%% file: paper1_figueroa-30jul.tex
\documentclass[12pt]{article}

\usepackage{graphicx,color,amsmath,amssymb,amsthm,rotating,stmaryrd,datetime,pifont,booktabs,enumitem}
\usepackage[mathscr]{eucal}

\setlist{nosep}

\newtheorem{theorem}{Theorem}

\newtheorem{lemma}[theorem]{Lemma}
\newtheorem{corollary}[theorem]{Corollary}

\newtheorem{remark}[theorem]{Remark}

\usepackage{tikz}

\renewcommand{\line}{{\tikz{\fill (0pt,0pt) circle (.75pt);\fill ( 0pt,2pt) circle (.75pt);\fill (0pt,-2pt) circle (.75pt)}}}
\newcommand{\plane}{{\tikz{\fill (0pt,0pt) circle (.8pt);\fill ( 1.5pt,2.5pt) circle (.8pt);\fill (3pt,0pt) circle (.8pt)}}}

\newcommand{\Label}{\label}

\newcommand\red[1]{{\color{red} #1}}

\definecolor{eggplant}{rgb}{0.68, 0.65, 0.32}

\newcommand{\XXX}{{$\overbrace{\phantom{XXXXXXXxxX}}^{{\mathcal D}}$}}
\newcommand{\YYY}{{$\overbrace{\phantom{XXXXXXXXXXXX}}^{\mathbb G'}$}}
\newcommand{\GGG}{$\left.\phantom{\begin{array}{c}X\\X\\X\\X\\X\\X\\\end{array}}\right\}\ \ {\mathbb G}\ \ \textup{(the fixed-\III-planes meeting }\gamma=\Bo{G}\textup{)}$}

\newcommand{\si}{\Sigma_\infty}
\newcommand{\li}{\ell_\infty}

\newcommand{\A}{\mathcal A}
\newcommand{\B}{\mathbf B}
\newcommand{\X}{\mathcal X}

\newcommand{\F}{\mathcal{F}}
\newcommand{\PG}{{\textup{PG}}}
\newcommand{\PGL}{{\textup{PGL}}}
\newcommand{\PGammaL}{\textup{P}\Gamma{\textup {L}}}
\newcommand{\Fqq}{\mathbb{F}_{q^2}}
\newcommand{\Fqqq}{\mathbb{F}_{q^3}}
\newcommand{\Fq}{\mathbb{F}_{q}}

\newcommand{\Q}{\mathcal{Q}}
\newcommand{\T}{\mathcal{T}}
\newcommand{\I}{\mathcal{I}}
\renewcommand{\H}{\mathbb{H}}
\renewcommand{\P}{\mathcal{P}}
\renewcommand{\S}{\mathbb{S}}
\renewcommand{\L}{\mathcal{L}}
\newcommand{\K}{\mathcal{K}}

\newcommand{\E}{\mathcal{E}}
\newcommand{\R}{\mathcal{R}}

\newcommand{\C}{\mathcal{C}}
\newcommand{\D}{\mathcal{D}}

\newcommand{\PI}{{\mathbb P}_{\mbox{\textup{\!\tiny $\I$}}}}
\newcommand{\PII}{{\mathbb P}_{\mbox{\textup{\!\tiny $\II$}}}}
\newcommand{\PIII}{{\mathbb P}_{\mbox{\textup{\!\tiny $\III$}}}}

\newcommand{\LI}{{\mathbb L}_{\mbox{\textup{\tiny $\I$}}}}
\newcommand{\LII}{{\mathbb L}_{\mbox{\textup{\tiny $\II$}}}}
\newcommand{\LIII}{{\mathbb L}_{\mbox{\textup{\tiny $\III$}}}}
\newcommand{\LFig}{{\mathbb L}_{\mbox{\textup{\tiny Fig}}}}

\newcommand{\qzero}{q\equiv0\pmod3}
\newcommand{\qone}{q\equiv1\pmod3}
\newcommand{\qtwo}{q\equiv-1\pmod3}

\newcommand{\piq}{\mathcal P_{2,q}}

\newcommand{\pifix}{\pi_{{\textup{\tiny{fix}}}}}

\newcommand{\Pifix}{\Pi_{{\textup{\tiny{Fix}}}}}
\newcommand{\SX}{\S_{\mbox{\textup{\tiny X}}}}

\newcommand{\SI}{\S_{\mbox{\textup{\tiny $\I$}}}}
\newcommand{\SII}{\S_{\mbox{\textup{\tiny $\II$}}}}
\newcommand{\SIII}{\S_{\mbox{\textup{\tiny $\III$}}}}
\newcommand\FIG{\mbox{FIG}}

\newcommand{\HI}{\H_{\mbox{\textup{\tiny $\I$}}}}
\newcommand{\HII}{\H_{\mbox{\textup{\tiny $\II$}}}}
\newcommand{\HIII}{\H_{\mbox{\textup{\tiny $\III$}}}}

\usepackage{stmaryrd}
\newcommand{\rbose }{\rrbracket}\newcommand{\lbose}{\llbracket}
\newcommand\Bo[1]{\mbox{$\lbose#1\rbose$}}

\makeatletter
\@ifpackageloaded{stix}{%
}{%
  \DeclareFontEncoding{LS2}{}{\noaccents@}
  \DeclareFontSubstitution{LS2}{stix}{m}{n}
  \DeclareSymbolFont{stix@largesymbols}{LS2}{stixex}{m}{n}
  \SetSymbolFont{stix@largesymbols}{bold}{LS2}{stixex}{b}{n}
  \DeclareMathDelimiter{\lDoubleBrace}{\mathopen} {stix@largesymbols}{"E8}%
                                            {stix@largesymbols}{"0E}
  \DeclareMathDelimiter{\rDoubleBrace}{\mathclose}{stix@largesymbols}{"E9}%
                                            {stix@largesymbols}{"0F}
}
\makeatother
\newcommand\figg[1]{\lDoubleBrace\hspace*{-0.12em} #1 \hspace*{-0.12em}\rDoubleBrace}

\newcommand\spr[1]{\langle #1,{#1}^{\sigma},{#1}^{\sigma^2}\rangle}

\newcommand{\g}{\mathsf g}
\newcommand{\m}{\mathsf n}
\newcommand{\III}{\textup{III}}
\newcommand{\II}{\textup{II}}
\renewcommand{\I}{\textup{I}}
\newcommand{\h}{{\texttt{\textup h}}}

\newcommand{\IIa}{{\textup{II}$_{\mbox{\textup{\scriptsize 1}}}$}}
\newcommand{\IIb}{{\textup{II}$_{\mbox{\textup{\scriptsize 2}}}$}}
\newcommand{\ha}{{{\texttt{\textup h$_{\mbox{\textup{\scriptsize 1}}}$}}}}
\newcommand{\hb}{{{\texttt{\textup h$_{\mbox{\textup{\scriptsize 2}}}$}}}}

\begin{document}

\title{The planar projectivity of $\PG(2,q^3)$ of order $3$ under field reduction}

\author{S.G. Barwick, Alice M.W. Hui and Wen-Ai Jackson
\date{}
}

\maketitle



\begin{abstract} Let $\phi $ be a collineation of $\PG(2,q^3)$ of order 3 which fixes  a plane of order $q$ pointwise. The points of   $\PG(2,q^3)$  can be partitioned into three types with respect to orbits of $\phi$: fixed points; points $P$ with $P,P^{\phi},P^{{\phi^2}}$ distinct and collinear; and points $P$ with $P,P^{\phi},P^{{\phi^2}}$ not collinear.  Under field reduction, the collineation $\phi$ corresponds to a projectivity  $\sigma$ of $\PG(8,q)$ of order 3. With respect to the field reduction and the orbits of $\sigma$, the points of   $\PG(8,q)$  can be partitioned into six types. This article looks at the projectivity  $\sigma$ in detail, and classifies and counts the fixed points, fixed lines and fixed planes.
The motivation is to give a description of the lines of the Figueroa projective plane in the $\PG(8,q)$ field reduction setting.
\end{abstract}

%

\section{Introduction}
\Label{001}

The Andr\'e/Bruck-Bose representation \cite{andre,BB1,BB2,segre} of a translation plane of order $q^n$ with kernel $\Fq$ in $\PG(2n,q)$ is a well-known use of field reduction. This representation has been used to study translation planes; including constructing new translation planes, classifying translation planes and studying structures in translation planes.
In this article we examine a non-translation plane in the field reduction setting, namely the Figueroa plane. We can use field reduction with the Figueroa plane  since we can identify the points of the Figueroa plane with the points of the Desarguesian plane $\PG(2,q^3)$.
We note that both translation and non-translation planes can be represented in a higher dimensional projective space using the  cone representation, which  is a generalization of the Andr\'e/Bruck-Bose representation, see Chapter 15 of \cite{JJB} and references therein. However the cone representation has not been widely used to study projective planes, whereas the field reduction representation has.
Our main result on Figueroa planes in Theorem~\ref{407} determines the representation of Figueroa lines in the $\PG(8,q)$ field reduction setting.

The Figueroa planes are of prime power order $q^3$, $q>2$, they were constructed in  \cite{figueroa} for the case $q\not\equiv1\pmod3$ and generalised to include the case $\qone$ in \cite{hering}.
Figueroa planes are one of the three infinite classes of finite projective planes that are not closely connected to translation planes; the other two classes are Hughes and Coulter-Matthews, see \cite[Thm 11.1/remark]{Johnson}.
An interesting open question is whether the construction can be generalised to construct a new class of non-translation planes.
The synthetic construction of the Figueroa plane given in  \cite{grundhofer} relies on a   collineation $\phi$ of $\PG(2,q^3)$ of order 3 which fixes pointwise an $\Fq$-plane.
This article studies the collineation $\phi$ in the $\PG(8,q)$ field reduction setting; the aim  is to create the framework we need to study the Figueroa plane under field reduction. Throughout this article we assume $q>2$.

Let $\mathbb S$ be a regular (Desarguesian) $2$-spread in $\PG(8,q)$, so $ \mathbb S$ is
a set of $q^6+q^3+1$ planes, called \emph{$\mathbb S$-planes}. Let $\mathbb H$ denote the set of $5$-spaces which contain $q^3+1$ planes of $\mathbb S$, then $\H$ contains $q^6+q^3+1$ $5$-spaces, called  \emph{$\mathbb H$-$5$-spaces}. Further, $\mathbb H$ forms a dual spread   (that is, each
$7$-space of $\PG(8,q)$ contains a unique $5$-space in $\mathbb H$).
As the $\mathbb S$-planes are pairwise disjoint, any two $\mathbb S$-planes span a $\mathbb H$-$5$-space, and  any two $\mathbb H$-$5$-spaces meet in
an $\mathbb S$-plane.
 Let $\mathcal I$ be the incidence structure with  \emph{points} the $q^6+q^3+1$
$\mathbb S$-planes;  \emph{lines} the  $5$-spaces of $\PG(8,q)$ that contain two (and so $q^3+1$) $\mathbb S$-planes; and  \emph{incidence} is inclusion.
Then $\mathcal I\cong\PG(2,q^3)$, and we call $\mathcal I$ the Bose representation of $\PG(2,q^3)$ in $\PG(8,q)$. This representation generalises the
 construction given by Bose \cite{bose} to represent $\PG(2,q^2)$ using a regular $1$-spread in $\PG(5,q)$.
The Bose representation is an example of the technique of field reduction. This idea goes back to Segre \cite{segre} who introduced regular spreads arising from field reduction. Field reduction  has been an area of much recent interest, see \cite{fieldreduction} for a survey.

  We use uppercase letters to denote points in $\PG(8,q)$   and uppercase
letters with a bar to denote points in $\PG(2,q^3)$.
 A point $\bar P$ of $\PG(2,q^3)$ corresponds in $\PG(8,q)$ to an $\S$-plane denoted $\Bo P$; and a line $\bar \ell$  of $\PG(2,q^3)$ corresponds to an $\mathbb H$-$5$-space denoted $\Bo\ell$.
   If $\bar \X$ is a set of points in $\PG(2,q^3)$, then the corresponding set of $\S$-planes is denoted $\Bo{\X}$, that is, $\lbose \X\rbose=\{\Bo{P}\,|\,\bar P\in\bar \X\}$. See \cite{BJW3} for a  detailed analysis of this representation.

  Let $\phi$ be a   collineation of $\PG(2,q^3)$ of order 3 which fixes pointwise an $\Fq$-plane which we denote by $\piq$.
We partition the points
 of $\PG(2,q^3)$ according to the orbits of $\phi $.
 A point $\bar P$ has Type $\I$, Type $\II$ and Type $\III$ respectively if the orbit of $\bar P$ under $\phi $ is a point, three collinear points, and three non-collinear points.
Dually,  a line $\bar \ell$ has Type $\I$, Type $\II$ and Type $\III$ respectively if the orbit of $\bar \ell$ under $\phi $ is a line, three concurrent lines, and three non-concurrent lines.
Thus, a point of $\piq$ has Type $\I$; the remaining points  on  lines of $\piq$ have Type $\II$, and the rest have Type $\III$. A line has Type $\I$ if it meets $\piq$ in $q+1$ points;  Type $\II$ if it meets $\piq$ in exactly one point, and  Type $\III$ if it is exterior to $\piq$.

In the Bose representation, we label the $\mathbb S$-planes and $\H$-$5$-spaces  of $\PG(8,q)$ with the same type as the corresponding points and lines of $\PG(2,q^3)$.
We use the notation: an $\SI$-, $\SII$-, $\SIII$-plane of $\PG(8,q)$ corresponds to a Type $\I$, $\II$, $\III$ point  of $\PG(2,q^3)$ respectively;
and an $\HI$-, $\HII$-, $\HIII$-$5$-space  of $\PG(8,q)$   corresponds to
a Type $\I$, $\II$, $\III$ line  of $\PG(2,q^3)$ respectively.

The collineation $\phi\in\PGammaL(3,q^3)$ induces a projectivity  $\sigma\in\PGL(9,q)$ of order 3 acting on $\PG(8,q)$.
If $\Pi $ is a subspace of $\PG(8,q)$, then $\sigma$ fixes the subspace $\spr{\Pi}$. By
 looking at the action of $\phi $ on the points of $\PG(2,q^3)$,
we have the following.
\begin{itemize}
\item If $\alpha$ is an $\SI$-plane, then
$\spr{\alpha}=\alpha$.
\item  If $\beta$ is an $\SII$-plane, then $\spr{\beta} $ is the unique $\HI$-$5$-space containing $\beta$.
\item  If $\gamma$ is an $\SIII$-plane,  then  $\spr{\gamma}=\PG(8,q)$.
\end{itemize}

As the $\mathbb S$-planes partition the points of $\PG(8,q)$, we   label points according to the $\S$-plane containing them.
That is, if a point $P$ of $\PG(8,q)$ lies in an $\SI$-plane (respectively $\SII$-plane or $\SIII$-plane), then we call $P$  a \emph{$\I$-point} (respectively \emph{$\II$-point} or \emph{$\III$-point}).
To help distinguish the setting, the phrase \emph{$X$-point} is reserved for points in $\PG(8,q)$, and the phrase \emph{Type-$X$-point} is reserved for points in $\PG(2,q^3)$.
 The number of $X$-points in $\PG(8,q)$ equals $q^2+q+1$ times the number of Type-$X$-points in $\PG(2,q^3)$, $X\in\{\I,\II,\III\}$. We further partition the points of $\PG(8,q)$
by looking at their orbit under $\sigma$. It is straightforward to verify that there are six possible types of points, and we use the following terminology for these.

\begin{itemize}
\item Let $A$ be a \I-point, then $\langle A,A^\sigma,A^{\sigma^2}\rangle$ is contained in a unique $\SI$-plane and either:
\begin{itemize}
\item
 $A$ is fixed by $\sigma$ and is called a  \emph{fixed-point};
 \item $A,A^\sigma,A^{\sigma^2}$ are distinct and collinear and $A$ is called a \emph{$\I_\line$-point}; or
 \item
 $A,A^\sigma,A^{\sigma^2}$ form a triangle and $A$ is called  a \emph{$\I_\plane$-point}.
 \end{itemize}

\item Let $B$ be a  \II-point, then $\langle B,B^\sigma,B^{\sigma^2}\rangle$ is contained in a unique $\HI$-$5$-space and either:
\begin{itemize}
\item

  $B,B^\sigma,B^{\sigma^2}$ are collinear  and $B$ is called
a  \emph{$\II_\line$-point}; or
\item  $B,B^\sigma,B^{\sigma^2}$  form a triangle  and $B$ is called a
 \emph{$\II_\plane$-point}.
 \end{itemize}
\item  Let $C$ be a   \III-point, then $C,C^\sigma,C^{\sigma^2}$  form a triangle, the plane $\langle C,C^\sigma,C^{\sigma^2}\rangle$ is not contained in any $\H$-$5$-space,  and $C$ is called a
 \emph{$\III_\plane$-point}.
 \end{itemize}

In this article, we study the projectivity  $\sigma$ in detail.  That is, we study the orbits of the projectivity $\phi$ after field reduction. This is a natural question to ask in the context of studying Figueroa planes in the setting of field reduction. Moreover this study is used to give  a  representation of the Figueroa lines in the field reduction setting  in Section~\ref{Sec6}.

The article proceeds as follows.
 In Section~\ref{003}, we classify the points, lines and planes of $\PG(8,q)$ fixed by $\sigma$. We count the fixed subspaces of each type, and look at the structure of the fixed subspaces. Writing $\g=\gcd(3,q-1)$, we show in Theorem~\ref{000} that there are exactly $\g(q^2+q+1)$ fixed points, and in Theorem~\ref{000} that these form $\g$ planes fixed pointwise by $\sigma$. Dually, there are $\g$ $5$-spaces fixed hyperplanewise by $\sigma$, and we study interesting properties of these hyperplanes. The importance of these hyperplanewise-fixed $5$-spaces is illustrated in Remark~\ref{rem9} and
 Theorem~\ref{333a} which show that they contain all the $\I_\line$-points and $\II_\line$-points; as well as all the fixed points and fixed lines when $q\not\equiv -1\pmod{3}$.
 In Theorem~\ref{294} we show that there are three types of fixed lines. The nice geometric structure of fixed lines is looked at in Section~\ref{sec-more},  in particular, Theorem~\ref{191} describes how they form linear congruences. We conclude this discussion on fixed subspaces by summarising the counting results in Tables~\ref{table16}, \ref{table14} and \ref{table18}.

In Section~\ref{Sec3}, we look at an application to linear sets of $\PG(2,q^3)$. In particular  each  $\sigma$-fixed subspace of $\PG(8,q)$ corresponds to a $\phi$-fixed linear set of $\PG(2,q^3)$.  We determine the structure of these linear sets.
In Section~\ref{Sec5}, we briefly consider the Bruck-Bose representation of $\PG(2,q^3)$ in $\PG(6,q)$. We first observe how the results in $\PG(8,q)$ have an interpretation in $\PG(6,q)$.
 We then describe an interesting difference that arises in the $\PG(6,q)$ model, namely that
the set of all \I-points and \II-points in $\PG(6,q)$ correspond precisely to the points of a \emph{single quadric} of $\PG(6,q)$. We show how this quadric gives an insight into the structure of the \II-points in $\PG(6,q)$.
Finally, Section~\ref{Sec6} looks at the Figueroa plane construction in the field reduction setting. We aim to find a representation for the lines of the Figueroa plane in the $\PG(8,q)$ setting. We show  in Theorem~\ref{407}  that when $q\not\equiv1\pmod3$, the Type $\III$ Figueroa lines can be represented in $\PG(8,q)$ using rational normal $3$-fold scrolls.

\section{Preliminaries}\label{sec-pre}

In this article, we will study the subspaces of $\PG(8,q)$  that are fixed by the collineation $\sigma$. For the remainder of the article, whenever we refer to a \emph{fixed subspace}, we mean a subspace of $\PG(8,q)$ that is fixed by the collineation $\sigma$.

\subsection{Notation}

We let $\Fq$ denote the finite field of order $q$ for $q$ a prime power. Further, $\Fq^*=\Fq\setminus\{0\}$.

The collineation $\sigma$ is a projectivity  of order 3, so points of $\PG(8,q)$ have orbits of size 1 or 3 under $\sigma$.
Suppose $\ell$ is a fixed line, then $\sigma$ induces a projectivity  of order 3 acting
on the points of $\ell$. Hence $\sigma$ fixes $x\in\{0,1,2,q+1\}$ points of $\ell$. As the non-fixed points of $\ell$ lie in orbits of size $3$, $q+1-x$ is divisible by $3$. Hence the number of fixed points on $\ell$ depends on the value of $q+1$ modulo $3$.
Throughout this article, we
 distinguish these three cases using the notation: $$q\equiv\m \pmod3,\quad\m \in\{-1,0,1\}.$$
Similarly, when looking at fixed planes, the value of $q^2+q+1$ modulo $3$ is important.
 If $\qtwo$ or $\qzero$, then $\gcd(3,q-1)=1$ and $q^2+q+1\equiv 1$.
  If $\qone$, then $\gcd(3,q-1)=3$ and $q^2+q+1\equiv 0$.
Throughout this article, we
 distinguish these two cases using the notation:   $$\g =\gcd(3,q-1),\quad\g \in\{1,3\}.$$

\subsection{The Bose representation of $\Fq$-lines and $\Fq$-planes of $\PG(2,q^3)$}\label{sec:bose}

 The representation of
$\Fq$-lines and $\Fq$-planes  of  $\PG(2,q^3)$ in the  Bose representation in $\PG(8,q)$ is known. For example, it is proved in    \cite{lunardon}
and also in \cite[Theorem 2.6]{fieldreduction} using field reduction techniques.
We need the follow details. Firstly,  an $\Fq$-line $\bar b$ of $\PG(2,q^3)$ corresponds in $\PG(8,q)$ to a set $\Bo{b}$ of $q+1$ $\S$-planes which form a 2-regulus $\T$. Note that a 2-regulus determines a unique  Segre variety ${\mathscr S}_{1;2}$, the ruling planes are the $q+1$  planes of $\T$, and the $q^2+q+1$ ruling lines meet every plane in $\T$.
Secondly, an $\Fq$-plane $\bar \pi$  of $\PG(2,q^3)$ corresponds in $\PG(8,q)$ to a set $\Bo{\pi}$ of $q^2+q+1$ $\S$-planes which form  one system of ruling planes of a   Segre variety  ${\mathscr S}_{2;2}$.
This Segre variety has a second system of
 $q^2+q+1$
 ruling planes. See \cite[Thm 25.5.5]{HT} for more details on the Segre variety ${\mathscr S}_{2;2}$. In particular note that two ruling planes in the same system are disjoint; and two ruling planes from different systems  meet in a unique point.

\section{Subspaces of $\PG(8,q)$ fixed by $\sigma$}\Label{003}

 In this section we study subspaces of $\PG(8,q)$ fixed by $\sigma$. Of particular interest are the fixed points, fixed lines, fixed planes and $5$-spaces that are fixed hyperplane-wise by $\sigma$.
For ease of reference, the counting results obtained in this section are summarised in Tables~\ref{table16}, \ref{table14} and \ref{table18}.

\subsection{Points fixed by $\sigma$}\Label{sec2.2}

The collineation $\phi $ acting on $\PG(2,q^3)$ fixes pointwise the $\Fq$-plane $\piq$. By Section~\ref{sec:bose},  $\piq$ corresponds  in $\PG(8,q)$ to the set $\Bo{\piq}$  which consists of the  $q^2+q+1$ $\SI$-planes. That is, the $\SI$-planes  form one system of ruling planes of a   Segre variety  ${\mathscr S}_{2;2}$. This Segre variety has a second system of
 $q^2+q+1$
 ruling planes which we call \emph{$\piq$-ruling-planes}.

The next result  uses coordinates to count the fixed points of $\sigma$ and to show that exactly $\g=\gcd(3,q-1)$ $\piq$-ruling planes are fixed.
A plane of $\PG(8,q)$ that is fixed pointwise by $\sigma$ is called a \emph{ptwise-fixed plane}.

\begin{theorem} \Label{000}
\begin{enumerate}
\item The number of  points of $\PG(8,q)$ fixed by $\sigma$ is $\g (q^2+q+1)$.
\item Each $\SI$-plane contains exactly $\g$ fixed points.
\item The fixed points   form exactly $\g$  ptwise-fixed planes, moreover these  are $\piq$-ruling planes.
\end{enumerate}
\end{theorem}

\begin{proof}
We use coordinates to prove this result and use the  notation: if $\bar P$ is a point in $\PG(2,q^3)$, then $\bar P$ has homogeneous coordinates $\vec{P}=(x,y,z)$ for some $x,y,z\in\Fqqq$, not all zero. Further,   $\vec{P}\equiv\lambda (x,y,z)$ for any $\lambda\in\Fqqq^*$.

By \cite[Theorem 2]{BJ}, any planar order 3 collineation $\phi$ of $\PG(2,q^3)$
is conjugate to $\phi_1$ or $\phi_1^2$, where $\phi_1\in\PGammaL(3,q^3)$ is:
\begin{equation}\label{eqn:0}
\begin{array}{rrcl}
\phi_1 \colon\  & (x,y,z)& \longmapsto& (z^q,x^q,y^q)\\
& [d,e,f] & \longmapsto & [f^q,d^q,e^q]
\end{array}.
\end{equation}
We will prove the result when
$\phi=\phi_1.$
 The proof when $\phi=\phi_1^2$ is similar.
The $\Fq$-plane fixed pointwise by $\phi $ is
$\piq=\{(x,x^q,x^{q^2})\,|\, x\in\Fqqq^*\}.$

Let $\tau$ be a primitive element of $\Fqqq$ over $\Fq$. So each element $x\in\Fqqq$ can be uniquely written as $x=x_0+x_1\tau+x_2\tau^2$ with $x_0,x_1,x_2\in\Fq$.
Define the following two maps
$$\begin{array}{rrclcrrcl}
\theta:&\Fqqq&\longrightarrow &\Fq^3&\phantom{XXXXXXXX}&\Theta:&\Fqqq^3&\longrightarrow &\Fq^9\\
&x&\mapsto&(x_0,x_1,x_2)&
&&(x,y,z)&\mapsto&(\theta(x),\theta(y),\theta(z)).
\end{array}
$$
If $\bar P$ is a point of $\PG(2,q^3)$ with homogeneous coordinates $\vec P=(x,y,z)$, then the corresponding $\mathbb S$-plane  in $\PG(8,q)$ contains the points with coordinates given by $$\Bo{P}=\{\Theta(\lambda x,\lambda y,\lambda z)\, |\, \lambda\in\Fqqq^*\}.$$

 We next determine the matrix for the collineation $\sigma \in\PGL(9,q)$. Let $\theta(\tau^q)=(a_0,a_1,a_2)$ and  $\theta({\tau^{2q}})=(b_0,b_1,b_2)$.
Let $A$ be the following $3\times3$ matrix over $\Fq$, and $M$ the following $9\times9$ matrix over $\Fq$ (where $O$ denotes the $3\times 3$ matrix of all zeroes).
\[
A=\begin{pmatrix}
1&0&0\\
a_0&a_1&a_2\\
b_0&b_1&b_2
\end{pmatrix},\quad\quad\quad M=\begin{pmatrix}
O&A&O\cr
O&O&A\cr
A&O&O
\end{pmatrix}.
\]
Observe that  for $x\in\Fqqq$, $\theta(x)A=\theta(x^q)$.
We have  $\sigma\colon\, (x_0,\ldots,x_8)\mapsto (x_0,\ldots,x_8)M$, since
\begin{eqnarray*}
\sigma(\Theta(x,y,z))&=&\sigma(\theta(x), \theta(y),\theta(z))\\
&=&(\,\theta(z)A, \,\theta(x)A,\,\theta(y) A\,)\\
&=&(\,\theta(z^q), \,\theta(x^q),\,\theta(y^q)\,)\\
&=&\Theta(z^q,x^q,y^q)
\end{eqnarray*} agreeing with $\phi=\phi_1$ which is given in (\ref{eqn:0}).

Recall that $\SI$-planes are fixed, whereas $\SII$-planes and $\SIII$-planes are not fixed and hence they are disjoint to their images.
Hence a fixed point must lie in an $\SI$-plane.
Let $\alpha$ be an $\SI$-plane, so for some given $x\in\Fqqq^*$, the points of $\alpha$ have homogeneous coordinates of form $\vec{Q_ t }=\Theta( t  x, t  x^q, t  x^{q^2})=(\,\theta( t  x),\,\theta( t  x^q),\,\theta( t  x^{q^2})\,),$ $ t \in\Fqqq^*$. The point  $Q_ t $
has image $\sigma(\vec Q_t)= (\, \theta( t ^q x),\,\theta(  t ^q x^q),\,\theta(  t ^q x^{q^2})\,)$.
 So $Q_t$ is a fixed point iff   there exists $\lambda\in\Fq$ with  $\theta( t ^qx^{q^i})=\lambda \theta( t  x^{q^i})$, $i=0,1,2$. As $\lambda
\in\Fq$, this holds iff  $ t ^qx^{q^i}=\lambda t  x^{q^i}$,  that is $ t ^q =\lambda t  $.  We rewrite this as $\theta( t )A=\lambda\theta( t )$, and so the point $Q_ t $ is fixed by $\sigma$ iff $\theta( t )$ is an eigenvector of $A$
with eigenvalue $\lambda$.

We compute
\[\begin{array}{rclcrclc}A^2&=& \begin{pmatrix} {\theta(1)}\\ {\theta(\tau^{q^2})}\\{\theta(\tau^{2q^2})}\end{pmatrix}&\quad \textup{and}\quad       A^3&=& \begin{pmatrix} {\theta(1)}\\ {\theta(\tau^{q^3})}\\{\theta(\tau^{2q^3})}\end{pmatrix}
         &=& \begin{pmatrix} {\theta(1)}\\ {\theta(\tau)}\\{\theta(\tau^{2})}\end{pmatrix}=I,
\end{array}
\]
hence $A$ has characteristic equation $\lambda^3-1$.
We show that the eigenvalue $\lambda=1$ has an eigenspace of dimension one. Now $\theta( t )=\theta( t ^q)$ iff  $ t = t ^q$. We can write $
t = t _0+ t _1\tau+ t _2\tau^2$, $ t _i\in\Fq$ and so $ t _1(\tau^q-\tau)+ t _2(\tau^{2q}-\tau^2)=(\tau^q-\tau)( t _1+ t _2(\tau^q+\tau))=0$
 from which we conclude that $ t _1= t _2=0$.  Thus the eigenspace has dimension 1 with basis $\{ (1,0,0)\}$. This corresponds to a fixed point in $\alpha$ of coordinates $\vec Q_1=\Theta(x,x^q,x^{q^2})$.

We now determine any remaining eigenvectors of $A$.
If $q\equiv 0\pmod 3$ then $\lambda^3-1=(\lambda-1)^3=0$, so $A$ has one distinct eigenvalue $\lambda=1$, with the eigenspace discussed above. Hence the $\SI$-plane $\alpha$ contains exactly one fixed point (and dually, exactly one fixed line).
If $q\not\equiv 0\pmod 3$, then $\lambda^3-1=(\lambda-1)(\lambda^2+\lambda+1)$.  The quadratic  $\lambda^2+\lambda+1$ is reducible with root $k$
iff $\{1,k,k^2\}$ forms a multiplicative subgroup of $\Fq^*$, that is, iff $3$ divides $q-1$, that is $q\equiv 1\pmod 3$.  In this case $A$ has three distinct eigenvalues. Hence $\alpha$ contains  exactly three fixed points, and dually three fixed lines, which necessarily form the vertices and sides of a triangle.  If $\qtwo$, then $\lambda^3-1$ has exactly one root in $\Fq$, so $\alpha$ contains  exactly one fixed point and one fixed line.
We conclude that there are $\g (q^2+q+1)$ fixed points, $\g $ in each $\SI$-plane, proving parts 1 and 2.

To prove part 3, let $\pi$ be a $\piq$-ruling plane, then each point of $\pi$ lies in a unique $\SI$-plane. As
 each $\SI$-plane is fixed by $\sigma$,  if $\pi$ is fixed by $\sigma$, then $\pi$  is fixed pointwise.
By part 1, $\sigma$ fixes exactly $\g (q^2+q+1)$ points, hence $\sigma$ fixes at most $\g $ $\piq$-ruling planes.
If $\g=1$, then the number of $\piq$-ruling planes is $q^2+q+1\equiv 1\pmod3$, so $\sigma$ fixes at least one
$\piq$-ruling plane.  Hence
 $\sigma$ fixes exactly one $\piq$-ruling plane, and this plane contains all the fixed points.
If $\g=3$, then  the number of $\piq$-ruling planes is  $q^2+q+1\equiv 0\pmod3$, and we conclude that $\sigma$ fixes either zero or three $\piq$-ruling planes.
Observe that $\sigma$ fixes
 pointwise the plane
 $
 \pi=\{(\theta(x),\theta(x^q),\theta(x^{q^2}))\bigm| x\in\Fqqq^*\}
$
of $\PG(8,q)$. Further, $\pi$ meets each $\S$-plane in $\Bo{\piq}$, so $\pi$ is  a $\piq$-ruling-plane. Hence $\sigma$ fixes exactly three $\piq$-ruling-planes, moreover these three planes contains all the fixed points.
\end{proof}

 \subsection{Lines fixed by $\sigma$}

In this section, we determine the fixed lines of $\sigma$. In Theorem~\ref{294}, we show that there are three types of fixed lines, which we call ptwise-fixed lines, fixed-\I-lines and fixed-\II-lines.
Further,
 in Lemma \ref{345}, we show that
the  ptwise-fixed lines are exactly the lines contained in a  ptwise-fixed-plane. In Lemma~\ref{345b}, we show that the fixed-\I-lines in a $\HI$-space form $\g$ 1-reguli.
We study the structure of
the fixed-\I-lines and fixed-\II-lines in more detail in Section \ref{sec-more}.

\begin{theorem}\Label{294}  Let $\ell$ be a fixed line of $\PG(8,q)$,  then $\ell$ is one of the following three types.
\begin{itemize}
\item $\ell$ is fixed pointwise, then $\ell$ meets $q+1$ $\SI$-planes and
is contained in a ptwise-fixed plane. We call this a
 \emph{ptwise-fixed line};

 \item $\ell$ contains   $\m +1$ fixed points and $q-\m $ $\I_\line$-points; $\ell$ is contained in an $\SI$-plane. We call this a
\emph{fixed-\I-line};

\item $\ell$ contains $\m +1$ fixed points and $q-\m $ $\II_\line$-points; $\ell$ meets $\m+1$ $\SI$-planes and $q-\m$  $\SII$-planes and is contained in a unique  $\H$-$5$-space which has Type $\I$. We call this a \emph{fixed-\II-line}.
\end{itemize}
\end{theorem}

\begin{proof}
Let $\ell$ be a fixed line of $\PG(8,q)$, then for a point $P\in\ell$, we have $P^\sigma,P^{\sigma^2}\in\ell$. Hence each point on $\ell$ is either a fixed point, a $\I_\line$-point or a $\II_\line$-point.

Suppose $\ell$ is fixed pointwise. Since the number of ptwise-fixed planes is less than $q+1$, $\ell$ lies on a ptwise-fixed plane, which is a $\piq$-ruling plane by Theorem~\ref{000}. So $\ell$ meets $q+1$ $\SI$-planes, which are contained in a unique $\HI$-$5$-space. Note that this $\HI$-$5$-space is the unique  $\H$-$5$-space containing $\ell$.

Suppose $\ell$ contains a $\I_\line$-point $P$, then $P,P^\sigma,P^{\sigma^2}\in\ell$ are distinct and lie in a common $\SI$-plane.
Thus $\ell$ lies in an $\SI$-plane, and so contains only fixed points and $\I_\line$-points.
As $\sigma$ is a projectivity  of $\PG(8,q)$, and $\ell$ is not fixed pointwise by $\sigma$,
  $\ell$ contains at most two fixed points. The remaining points of $\ell$ are partitioned  into orbits of size 3 under $\sigma$.
If $q\equiv\m \pmod 3$, $\m \in\{-1,0,1\}$, then the number of points of $\ell$ modulo 3 is $\m +1\in\{0,1,2\}$, so $\ell$ contains exactly $\m +1$ fixed-points, and $q-\m$ $\I_\line$-points.

Suppose $\ell$ contains a $\II_\line$-point, then by the above argument, $\ell$ does not contain a $\I_\line$-point.  That is, $\ell$ contains only fixed points and $\II_\line$-points. A similar argument to the previous case shows that $\ell$ contains exactly $\m +1$ fixed-points, and $q-\m$ $\II_\line$-points. Let $Q\in\ell$ be a $\II_\line$-point. Then $Q$ lies in a unique $\SII$-plane $\beta$ which lies in a unique $\HI$-$5$-space $\Sigma_5$. We have $Q^\sigma\in\beta^\sigma\subset\Sigma_5$ and $Q^{\sigma^2}\in\beta^{\sigma^2}\subset\Sigma_5$. Hence $\ell$ lies in $\Sigma_5$, and $\ell$ does not lie in any other $\H$-$5$-space.
  \end{proof}

We can now describe the points in $\SI$-planes.

  \begin{corollary}\Label{294b}
  \begin{enumerate}
  \item Each $\SI$-plane contains exactly
\begin{itemize}
\item $\g $ fixed points, $\g (q-\m )$ $\I_\line$-points and $q^2+q+1-\g (q-\m +1)$ $\I_\plane$-points;
\item $\g $ fixed lines which are all fixed-\I-lines.
\end{itemize}
\item \label{item3249b}  There are exactly $\g (q-\m )(q^2+q+1)$ $\I_\line$-points and $(q^2+q+1-\g (q-\m +1))(q^2+q+1)$ $\I_\plane$-points in $\PG(8,q)$.

 \end{enumerate}
\end{corollary}

\begin{proof}
Let $\alpha$ be an $\SI$-plane. By Theorem~\ref{000}, $\alpha$ contains exactly $\g$ fixed points. As $\alpha$ is fixed,
$\sigma$ induces a collineation acting on $\alpha$, this collineation has $\g$ fixed points, and so has $\g$ fixed lines. By Theorem~\ref{294}, these fixed lines are fixed-\I-lines and part 1 follows. Part 2 follows as there are $q^2+q+1$ $\SI$-planes.
\end{proof}

 We next count the number of ptwise-fixed lines and fixed-\I-lines in $\PG(8,q)$. The  fixed-\II-lines are counted later in Corollary~\ref{333c}.

\begin{lemma}\Label{345}
\begin{enumerate}
\item The number of ptwise-fixed lines is $\g (q^2+q+1)$, they comprise the lines of the $\g $ ptwise-fixed planes.
\item
 The   number of fixed-\I-lines is  $\g (q^2+q+1)$, consisting of $\g $ in each $\SI$-plane.
 \end{enumerate}
\end{lemma}

\begin{proof}
By Theorem~\ref{000} and Theorem~\ref{000}, there are exactly $\g (q^2+q+1)$ fixed points which form $\g $ ptwise-fixed planes. Hence there are $\g (q^2+q+1)$ ptwise-fixed lines, proving part 1.
By Theorem~\ref{294}, each fixed-\I-line is contained in an $\SI$-plane.
By Corollary~\ref{294b},  every $\SI$-plane contains exactly $\g $  fixed-\I-lines. As there are $q^2+q+1$ $\SI$-planes, there are exactly $\g (q^2+q+1)$ fixed-\I-lines.
 \end{proof}

 \begin{lemma}\Label{345b} Let $\Sigma_5$ be an $\HI$-$5$-space, then $\Sigma_5$ contains exactly $\g$ ptwise-fixed lines and $\g(q+1)$ fixed-\I-lines. Further,   the fixed-\I-lines contained in $\Sigma_5$ form $\g $ 1-reguli, each of whose opposite 1-regulus contains exactly $\m +1$ ptwise-fixed lines. If $\g =3$, then the three 1-reguli pairwise meet in a ptwise-fixed line.
 \end{lemma}

 \begin{proof}
Let $\Sigma_5$ be an $\HI$-$5$-space. So $\Sigma_5$ corresponds in $\PG(2,q^3)$ to a Type-\I-line which contains exactly $q+1$ Type-\I-points forming an $\Fq$-line.    Hence there are exactly $q+1$ $\SI$-planes in $\Sigma_5$, and they form a 2-regulus $\T$.  By Theorem~\ref{000}, $\Sigma_5$ contains $\g (q+1)$ fixed points; these lie $\g$ in each $\SI$-plane and $q+1$ on each of the $\g$ ptwise-fixed planes.  Hence each ptwise-fixed plane meets $\Sigma_5$ in a line, that is the $\g (q+1)$ fixed points in $\Sigma_5$ form $\g $ ptwise-fixed lines. Moreover,
these $\g$ lines meet each of the $\SI$-planes in $\Sigma_5$ in a point, so they are ruling lines of the 2-regulus $\T$.

By Corollary~\ref{294b}, there are $\g (q+1)$ fixed-\I-lines in $\Sigma_5$, we  show that these form $\g$ $1$-reguli.  The map $\sigma$ induces a collineation acting on the fixed subspace $\Sigma_5$ that fixes $\g$ lines of $\Sigma_5$ pointwise, and so dually, fixes $\g$ 3-spaces of $\Sigma_5$ `4-space-wise'.
First suppose $\g =1$, so there is a unique 3-space $\Pi_3$ in $\Sigma_5$ such that all 4-spaces of $\Sigma_5$ containing $\Pi_3$ are fixed.
 Let $\alpha$ be an $\SI$-plane contained in $\Sigma_5$. As both $\alpha$
and $\Pi_3$ are fixed, the  subspace $\alpha\cap\Pi_3$ is   fixed. It follows from the dimension  theorem  that $\alpha\cap\Pi_3$ is non-empty.
  Suppose  $\alpha\cap\Pi_3$ is a point, denoted $A$. Let $\ell_0,\ldots,\ell_q$ be the $q+1$ lines of $\alpha$ through $A$. As $\Pi_3$ is fixed `4-space-wise' in $\Sigma_5$,  the 4-space $\langle \ell_i,\Pi_3\rangle$, $i=0,\ldots,q$, is fixed, and so  meets the fixed plane $\alpha$ in a fixed line. That is, $\alpha$ contains $q+1$ fixed lines. As $q+1>1=\g$, this is a contradiction to Corollary~\ref{294b}. Hence $\dim(\alpha\cap\Pi_3)\geq1$, so each of the $q+1$ $\SI$-planes in $\Sigma_5$ is either contained in $\Pi_3$ or meets $\Pi_3$ in a line.
Suppose $\Pi_3$ contains an $\SI$-plane $\alpha$. Let $\alpha'$ be another $\SI$-plane;  the above argument shows that there exists a line  $\ell$ contained in $\alpha'\cap\Pi_3$. As $\alpha\cap\alpha'=\emptyset$, we have dim$\langle\alpha,\ell\rangle=4$, a contradiction. Thus $\Pi_3$ does not contain an $\SI$-plane.
 Hence $\Pi_3$ meets each $\SI$-plane in a fixed line which is a  fixed-\I-line.
 By Corollary~\ref{294b}, there are exactly $q+1$ fixed-\I-lines in $\Sigma_5$, one in each $\SI$-plane.
Thus, the $q+1$ fixed-\I-lines in $\Sigma_5$ lie in the intersection of the $3$-space $\Pi_3$ and the 2-regulus
$\T$, hence they form a 1-regulus. By Theorem~\ref{294}, if $\qzero$, then the unique ptwise-fixed line in $\Pi_3$ meets each fixed-\I-line and so
is a line of the opposite 1-regulus; and if $\qtwo$, then the unique ptwise-fixed line in $\Pi_3$ is disjoint from the 1-regulus.

If $\g =3$, then
a similar argument shows that the fixed-\I-lines in $\Sigma_5$ form three
1-reguli. Moreover, each of the opposite reguli contains two of the three
ptwise-fixed lines contained in $\Sigma_5$.
 \end{proof}

 \subsection{$5$-spaces fixed hyperplane-wise by $\sigma$}

By Theorem~\ref{000}, $\sigma$ fixes pointwise exactly $\g $ planes
of $\PG(8,q)$. Dually,  $\sigma$  fixes exactly $\g $ $5$-spaces hyperplane-wise, we call such a $5$-space a \emph{hwise-fixed-$5$-space}.  Note that a hwise-fixed-$5$-space is not an $\H$-$5$-space.
The structure of the fixed subspaces contained in a hwise-fixed-$5$-space is both interesting in its own right, as well as being an important step towards classifying fixed lines and planes of $\PG(8,q)$.
In this section we show that every $\I_\line$-point, $\II_\line$-point, fixed-\I-line and fixed-\II-line is contained in a hwise-fixed-$5$-space.

\begin{lemma} \Label{900}
A hwise-fixed-$5$-space   meets
each $\SI$-plane  in  exactly a fixed-\I-line, meets each $\SII$-plane in exactly a $\II_\line$-point, and is disjoint from each $\SIII$-plane.
\end{lemma}

\begin{proof} Let $\Pifix$ be a hwise-fixed-$5$-space. The set of all $\SI$-planes spans $\PG(8,q)$, so there is an $\SI$-plane $\alpha$ that is not contained in $\Pifix$.
 By Corollary \ref{294b}, $\alpha$ contains exactly $\g$ fixed points and $\g$ fixed lines which are fixed-\I-lines. Further, by Theorem~\ref{294}, points in $\alpha$ not on a fixed line are $\I_\plane$-points. Hence  there is a $\I_\plane$-point $A\in\alpha$ with $A\notin\Pifix$.
As $\Pifix$ is hwise-fixed, the $6$-space $\Sigma_6=\langle \Pifix,A\rangle$ is fixed. Hence  $A^\sigma,A^{\sigma^2}\in\Sigma_6$, and so $\alpha=\langle A,A^\sigma,A^{\sigma^2}\rangle\subset \Sigma_6$. As $\alpha$ is not contained in $\Pifix$,  $\alpha\cap\Pifix$ is a   line. As both $\alpha$ and $\Pifix$ are fixed, the line $\alpha\cap\Pifix$ is fixed. By Corollary~\ref{294b}, $\alpha\cap\Pifix$ is a fixed-\I-line.
We conclude that every $\SI$-plane meets $\Pifix$ in at least a line,  and so  $\Pifix$ contains at least $r=(q^2+q+1)(q+1)$ \I-points.

Let $\beta$ be an $\SII$-plane and let $\Sigma_5$ be the unique $\HI$-$5$-space containing $\beta$. As $\Sigma_5$ contains $q+1$ $\SI$-planes and $\Pifix$ meets each $\SI$-plane in at least a line,
$\Sigma_5\cap\Pifix$ has dimension at least 3. By  the dimension theorem, each 3-space of $\Sigma_5$ meets every $\S$-plane of $\Sigma_5$.
Hence $\beta$ meets $\Pifix$ in at least a point.
Thus $\Pifix$ contains at least $s=(q^2+q+1)(q^3-q)$ \II-points. As $r+s$ equals the total number of points in a $5$-space, $\Pifix$ meets each $\SI$-plane  in exactly a fixed-\I-line, meets each $\SII$-plane in  exactly  one point, and is disjoint from each $\SIII$-plane. Moreover, we can
conclude that each $\HI$-$5$-space meets $\Pifix$ in a 3-space.

It remains to show that every \II-point in $\Pifix$ is a $\II_\line$-point.
 Let $B$ be a \II-point in $\Pifix$, so $B$ lies in a unique $\SII$-plane
and so in a unique  $\HI$-$5$-space denoted $\Sigma_5$. We suppose $B$ is
a $\II_\plane$-point and work to a contradiction. So
 $\langle B,B^\sigma,B^{\sigma^2}\rangle $ is a fixed plane which  lies in $\Sigma_3=\Sigma_5\cap\Pifix$. By the preceding paragraph, $\Sigma_3$ is a 3-space containing $q+1$ fixed-\I-lines. By Lemma~\ref{345b},   the $q+1$ fixed-\I-lines in $\Sigma_3$  form a 1-regulus $\R$, whose opposite regulus $\R'$ contains $\m +1$ ptwise-fixed lines.

 Let $ \X=\langle B,B^\sigma,B^{\sigma^2}\rangle \cap\R$, then the points in $ \X$ form either  a non-degenerate conic or two lines.  If $\X$ is a non-degenerate conic $\C$, then $\C$ has $q+1$ points, so $\C$ has exactly one point on each line of $\R$ and one point on each line of $\R'$. Since $\C=\{\ell \cap \langle B,B^\sigma,B^{\sigma^2}\rangle \,|\,\ell \in \R\}$ and each line $\ell\in\R$ is fixed, $\C$ is pointwise fixed, and so every line of $\R'$ is fixed. This   contradicts Lemma~\ref{345b} which shows $\R'$ contains exactly  $\m+1$ fixed lines.  Hence $\X$ is not a non-degenerate conic, so $\X$ consists of two lines, $\ell,m$ with $\ell\in\R$, $m\in\R'$. As each line in $\R$ is fixed, the line $m$ is fixed pointwise. If $\qtwo$, this contradicts
Lemma~\ref{345b} which shows that $\R'$ has no pointwise fixed lines.
If $q\not\equiv -1\pmod3$, then $\sigma$ induces a central collineation acting on the fixed plane $\langle B,B^\sigma,B^{\sigma^2}\rangle $ with axis $m$ and centre a point $V$. Thus $V,B,B^\sigma,B^{\sigma^2}$ are collinear, a contradiction. Hence every \II-point in $\Pifix$ is a $\II_\line$-point.
 \end{proof}

\begin{theorem} \Label{333} Let $\Pifix$ be a hwise-fixed-$5$-space. Then
$\Pifix$ contains exactly
\begin{itemize}
\item  $(\m +1)(q^2+q+1)$ fixed-points, $(q-\m )(q^2+q+1)$ $\I_\line$-points and $(q^3-q)(q^2+q+1)$ $\II_\line$-points.
\item $(\m +1)(q^2+q+1)$ ptwise-fixed lines, $q^2+q+1$ fixed-\I-lines and
$(q^3-q)(q^2+q+1)/(q-\m )$ fixed-\II-lines.
\item $\m +1$ ptwise-fixed planes, and is disjoint from the remaining $\g
-(\m +1)$ ptwise-fixed planes.
\end{itemize}
\end{theorem}

\begin{proof} There are $q^2+q+1$ $\SI$-planes and $(q^3-q)(q^2+q+1)$ $\SII$-planes. So by Lemma~\ref{900}, $\Pifix$ contains $q^2+q+1$ fixed-\I-lines and $(q^3-q)(q^2+q+1)$ $\II_\line$-points. By Theorem~\ref{294}, each fixed-\I-line contains $\m +1$ fixed points and $q-\m $ $\I_\line$-points, hence $\Pifix$ contains exactly
$(\m +1)(q^2+q+1)$ fixed points and
$(q-\m )(q^2+q+1)$ $\I_\line$-points. This accounts for all the points in
$\Pifix$.

 By Theorem~\ref{000}, the fixed points of $\sigma$ are exactly the
points on the $\g $ ptwise-fixed planes. As $\Pifix$ contains $(\m +1)(q^2+q+1)$ fixed points, $\Pifix$ contains $\m +1$ ptwise-fixed planes and is disjoint from the remaining $\g -(\m +1)$ ptwise-fixed planes. Similarly, the
ptwise-fixed lines are precisely the lines in a ptwise-fixed plane, so $\Pifix$ contains $(\m +1)(q^2+q+1)$ ptwise-fixed lines.

  If $B$ is a $\II_\line$-point in $\Pifix$, then the fixed-\II-line $\langle B,B^\sigma,B^{\sigma^2}\rangle$ lies in $\Pifix$.
Further, if two fixed-\II-lines meet, then they meet in a fixed-point. By
Theorem~\ref{294}, a fixed-\II-line contains $q-\m$ $\II_\line$-points. As there are $(q^3-q)(q^2+q+1)$  $\II_\line$-points in $\Pifix$,  the number of fixed-\II-lines in $\Pifix$ is $(q^3-q)(q^2+q+1)/(q-\m )$.
\end{proof}

\begin{remark}\label{rem9}\textup{
Note  that Theorem \ref{333} says that if $\qtwo$, then the unique ptwise-fixed plane and the unique hwise-fixed-$5$-space are disjoint. If $\qzero$, then  the unique ptwise-fixed plane is contained in the unique hwise-fixed-$5$-space. If $\qone$, then each of the three hwise-fixed-$5$-spaces is the span of two of the three ptwise-fixed planes; and we next show that they pairwise meet in exactly a ptwise-fixed plane.}
\end{remark}

\begin{corollary}\Label{334}
If $\g =3$, then
\begin{enumerate}
\item the three hwise-fixed-$5$-spaces have no common point,
\item any two hwise-fixed-$5$-spaces meet in exactly a ptwise-fixed plane.
\end{enumerate}
\end{corollary}

\begin{proof} Suppose $\g=3$ and
let $\Pifix$, $\Pifix'$ be two hwise-fixed-$5$-spaces. By Theorem~\ref{333}, each contains two of the three ptwise-fixed planes and is disjoint from the third. By Theorem~\ref{000}, the three ptwise-fixed-planes are $\piq$-ruling planes and so are pairwise skew. Hence $\Pifix\cap\Pifix'$ contains a ptwise-fixed plane, denote this by $\pifix$; and $\Pifix\cap\Pifix'$ contains no other fixed points. Suppose $\Pifix\cap\Pifix'$ contains a further point $P\notin\pifix$. As $P$ is not fixed, by Theorem~\ref{333}, $P$ is  an $X_\line$-point for some $X\in\{\I,\II\}$. The $3$-space $\langle P,\pifix\rangle$  contains the fixed-$X$-line $PP^\sigma$. By Theorem~\ref{294}, the line $PP^\sigma$ contains two fixed points, contradicting $\Pifix\cap\Pifix'$ containing no fixed points outside $\pifix$. Hence $\Pifix\cap\Pifix'$ is exactly a ptwise-fixed plane.

\end{proof}

\begin{theorem} \Label{333a}
\begin{enumerate}

\item Every $\I_\line$-point and every $\II_\line$-point   is contained in a unique hwise-fixed-$5$-space.
\item Every fixed-\I-line and fixed-\II-line is contained in a unique hwise-fixed-$5$-space.
\end{enumerate}
\end{theorem}

\begin{proof}
Simple counting using Corollary~\ref{294b}(\ref{item3249b}), Theorem~\ref{333} and Corollary~\ref{334} shows that every $\I_\line$-point lies in a unique hwise-fixed-$5$-space.

Let $B$ be a $\II_\line$-point, so  $\ell_B=\langle B,B^\sigma,B^{\sigma^2}\rangle $ is a  fixed-\II-line.
We will show that $B$ lies in a hwise-fixed-$5$-space $\Pifix $. We look at the three cases $\m=0,1,2$ separately.

Suppose
  $\m =0$, then  $\g =1$ and there is a unique hwise-fixed-$5$-space denoted $\Pifix$.
By Theorem~\ref{294}, the fixed-\II-line $\ell_B$ contains a unique fixed-point, denoted  $P$.  By Theorem~\ref{000}, the fixed point $P$ lies in an $\SI$-plane, which contains a unique fixed-\I-line denoted $\ell$. By Theorem~\ref{294} and Corollary~\ref{294b},  $P\in\ell$. By the above paragraph, $\ell$ is contained in $\Pifix$.
  Let $\beta=\langle  \ell_B,\ell\rangle$, as $\beta$ contains  two fixed lines  $\ell$ and $\ell_B$, $\beta$ is fixed by $\sigma$. Hence $\sigma$ induces a collineation on $\beta$ with at least two fixed lines, so dually $\beta $ contains at least two fixed points. Let $Q\neq P$ be a second fixed point in $\beta$. By Theorem~\ref{294}, as $\m=0$, the line $PQ$ is a ptwise-fixed line,
  so $\ell\neq PQ$, and we have $\beta=\langle Q,\ell\rangle$.
By
Theorem~\ref{000}, there are exactly $q^2+q+1$ fixed points, which by Theorem~\ref{333}, all lie in $\Pifix$. Hence $Q\in\Pifix$ and so
 $\beta=\langle \ell_B,\ell \rangle=\langle Q,\ell\rangle$  is contained
in $\Pifix$. That is,  $B
 \in\Pifix$ as required.

Suppose $\m =-1$,  then  $\g =1$ and there is a unique hwise-fixed-$5$-space denoted $\Pifix$. Suppose the point $B$   does not lie in the hwise-fixed-$5$-space $\Pifix$. As $\Pifix$ is fixed hyperplane-wise, the $6$-space $\langle B,\Pifix\rangle$ is fixed, and so contains the line $\ell_B$. Hence   $\ell_B\cap \Pifix$  is a fixed point, a contradiction  as $\Pifix $ has no fixed points by Theorem~\ref{333}. Hence $B
 \in\Pifix$ as required.

Suppose $\m=1$, so $\g =3$, and denote the three hwise-fixed-$5$-spaces by
$\Pi_{{\textup{\tiny{Fix},k}}}$, $k=0,1,2$.
Suppose $B$   does not lie in any hwise-fixed-$5$-space. Then  $\langle B,\Pi_{{\textup{\tiny{Fix},k}}}\rangle$,  $k=0,1,2$,   is a fixed 6-space, so  contains the line $\ell_B$. Hence    $P_k=\ell_B\cap \Pi_{{\textup{\tiny{Fix},k}}}$   is a fixed point, $k=0,1,2$.
By Corollary~\ref{334},
$\Pi_{{\textup{\tiny{Fix},0}}}\cap\Pi_{{\textup{\tiny{Fix},1}}}\cap\Pi_{{\textup{\tiny{Fix},2}}}=\emptyset$. Hence $|\{P_0,P_1,P_2\}|\geq2$, so $\ell_B=P_iP_j$ for some $i\neq j$. By Theorem~\ref{294}, $\ell_B$ contains exactly two fixed points, so by Theorem~\ref{000}, $\ell_B$ meets two ptwise-fixed planes.  By Theorem~\ref{333}, each hwise-fixed-$5$-space contains two ptwise-fixed planes, so  the line $\ell_B=P_iP_j$ is  contained in one of the  hwise-fixed-$5$-spaces, a contradiction. Hence $B$ lies in a hwise-fixed-$5$-space as required. By Corollary \ref{334}, this hwise-fixed-$5$-space is unique.
 This completes the proof of part 1. Part 2  follows from part 1.
  \end{proof}

Theorem~\ref{333a} illustrates the importance of the hwise-fixed-$5$-space(s). Note that by Remark~\ref{rem9}, if $q\not\equiv1\pmod{3}$, then all fixed points and all fixed lines are contained in a hwise-fixed-$5$-space.

We  conclude this section by describing the points in $\SII$- and $\SIII$-planes.

\begin{corollary}\Label{333b1}
\begin{enumerate}
\item Each $\SII$-plane contains
 exactly  $\g $ $\II_\line$-points and $q^2+q+1-\g $ $\II_\plane$-points;
and no fixed lines.
\item There are exactly $\g(q^2+q+1)(q^3-q)$ $\II_\line$-points and $(q^2+q+1-\g)(q^2+q+1)(q^3-q) $ $\II_\plane$-points in $\PG(8,q)$.
\item Each $\SIII$-plane contains
 exactly   $q^2+q+1$ $\III_\plane$-points;
and no fixed lines.
\item There are exactly $(q^2+q+1)q^3(q-1)^2(q+1)$   $\III_\plane$-points  in $\PG(8,q)$.
\end{enumerate}
\end{corollary}

\begin{proof} Let $\beta$ be an $\SII$-plane, by Lemma~\ref{900}, each of
the $\g$ hwise-fixed-$5$-spaces meets $\beta$ in a $\II_\line$-point. By Theorem~\ref{333a}, every $\II_\line$-point lies in   a unique  hwise-fixed-$5$-space.
Hence $\beta$ contains exactly $\g $ $\II_\line$-points, the remaining points are $\II_\plane$-points. By Theorem~\ref{294},  $\beta$ does not contain any fixed lines. This completes part 1.
Part 2 follows as there are $(q^2+q+1)(q^3-q)$ Type $\II$ points in $\PG(2,q^3)$, so $\SII$-planes in $\PG(8,q)$. It  follows from Theorem~\ref{294} that an $\SIII$-plane does not contain a fixed line. Hence $\SIII$-planes contain only $\III$-points, proving part 3. Part 4 follows as there are $q^3(q-1)^2(q+1)$  Type $\III$ points in $\PG(2,q^3)$, so $\SIII$-planes in $\PG(8,q)$.
\end{proof}

\
\subsection{Fixed lines and linear congruences}\Label{sec-more}

Theorem~\ref{333a} is a powerful result, and allows us to count the fixed-\II-lines (Corollary~\ref{333c}) and to
classify the  fixed lines in each $\H$-$5$-space  (Corollary~\ref{333e}). Moreover we use it to classify the points and fixed lines in the intersection of a hwise-fixed-$5$-space and a $\H$-$5$-space (Corollary~\ref{221}),
 and show that
a linear congruence naturally arises from these fixed lines (Theorem~\ref{191}).

\begin{corollary}\Label{333c}
 The number of fixed-\II-lines is  exactly  $\g (q^3-q)(q^2+q+1)/(q-\m )$.
\end{corollary}
\begin{proof}
By Theorem~\ref{333a}, every fixed-\II-line is contained in a  unique hwise-fixed-$5$-space. The number of fixed-\II-lines in a hwise-fixed-$5$-space is given in Theorem~\ref{333}. As there are $\g$ hwise-fixed-$5$-spaces, the result follows.
\end{proof}

\begin{corollary}\Label{333e}
\begin{enumerate}
\item  An  $\HI$-$5$-space contains  exactly $\g$ ptwise-fixed lines, $\g(q+1)$ fixed-\I-lines (which form $\g$ 1-reguli) and  $\g(q^3-q)/(q-\m)$ fixed-\II-lines.
\item  An  $\HII$-$5$-space contains  exactly $\g$ fixed lines which are all fixed-\I-lines.
\item  An  $\HIII$-$5$-space contains no fixed lines.
\end{enumerate}
\end{corollary}

\begin{proof} Let  $\Sigma_5$ be a $\HI$-$5$-space, so $\Sigma_5$ contains $q+1$ $\SI$-planes
and $q^3-q$ $\SII$-planes.
By Corollary~\ref{333b1}, the  number of $\II_\line$-points in $\Sigma_5$ is $\g(q^3-q)$. Each lies on a unique fixed-\II-line which necessarily lies in $\Sigma_5$. By Theorem~\ref{294}, a fixed-\II-line contains $q-\m$ $\II_\line$-points, hence there
are $\g(q^3-q)/(q-\m)$ fixed-\II-lines in $\Sigma_5$.  Lemma~\ref{345b} proves the remainder of part 1.

By Theorem~\ref{294}, each fixed-\I-line  lies in an $\SI$-plane. By Lemma~\ref{345}, each $\SI$-plane contains exactly $\g$ fixed-\I-lines.
Hence if $\Sigma_5$ is  an  $\H$-$5$-space, then  the number of fixed-\I-lines in $\Sigma_5$ is $\g$ times the number of $\SI$-planes in $\Sigma_5$.
 By Theorem~\ref{294},
$\HII$- and $\HIII$-$5$-spaces contain no ptwise-fixed line or fixed-\II-lines,
as any ptwise-fixed plane is contained in a unique $\HI$-$5$-space but not any $\HII$- and $\HIII$-$5$-spaces. Parts 2 and 3 follow.
\end{proof}

\begin{corollary}\Label{221}
\begin{enumerate}
\item A hwise-fixed-$5$-space and an $\HI$-$5$-space meet in a $3$-space which contains  exactly
\begin{itemize}
\item $(q+1)(\m +1)$ fixed points, $(q+1)(q-\m )$ $\I_\line$-points and $q^3-q$  $\II_\line$-points
\item $q+\m+2+(q^3-q)/(q-\m )$ fixed lines, consisting of $\m +1$ ptwise-fixed lines, $q+1$ fixed-\I-lines which form a 1-regulus, and $(q^3-q)/(q-\m )$ fixed-\II-lines.
\end{itemize}
 \item  A hwise-fixed-$5$-space and an $\HII$-$5$-space meet in  a plane which contains exactly
 \begin{itemize}
\item $\m +1$ fixed points, $q-\m $ $\I_\line$-points  and $q^2$ $\II_\line$-points;
 \item $1$ fixed line, namely a fixed-\I-line.
\end{itemize}
  \item  A hwise-fixed-$5$-space and an $\HIII$-$5$-space meet in  a plane
which contains only  $\II_\line$-points and contains no fixed lines.
  \end{enumerate}
\end{corollary}

\begin{proof} The numbers of points  and fixed lines of each type contained in an $\S$-plane  are computed in Corollaries~\ref{294b} and \ref{333b1}. The number of fixed lines of each type contained in an $\H$-$5$-space is computed in Corollary~\ref{333e}. It is then straightforward to compute the numbers of points   of each type contained in an $\H$-$5$-space.
By Theorem~\ref{333}, of the $\g$ fixed points contained in each $\SI$-plane, exactly $\m+1$ lie in a given hwise-fixed-$5$-space.
By
Theorem~\ref{333a}, the $\I_\line$-points, $\II_\line$-points, fixed-\I-lines and fixed-\II-lines all lie in some hwise-fixed-$5$-space. Moreover, a hwise-fixed-$5$-space contains no $X_\plane$-points, $X\in\{\I,\II,\III\}$. This gives the numbers stated in the result.
\end{proof}

By Corollary~\ref{221}, an $\HI$-$5$-space meets a hwise-fixed-$5$-space in a $3$-space that contains exactly $\m+1$ ptwise-fixed lines.
In the following theorem, we show that the fixed-\I-lines and fixed-\II-lines in such a $3$-space are linear congruences, which are defined as a set of lines whose line coordinates (Pl\"ucker coordinates) satisfy two linear independent equations. See \cite[Sections 15.2, 15.4 \& Table 15.10]{hirsch2} for more details on linear congruences.
In particular, we have the following geometric descriptions  in $\PG(3,q)$.  An \emph{elliptic linear congruence} is a regular line spread.
 A \emph{parabolic linear congruence} is a set of $q^2+q+1$ lines consisting of an axis $\ell$ and a set $\L$ of $q^2+q$ lines  such that: each line in $\L$ meets $\ell$; any two distinct lines in $\L$ are either disjoint or meet in a point of $\ell$; and $\L$ contains a regulus. 
A \emph{hyperbolic linear congruence} is the set of $(q+1)^2$ lines meeting a pair of disjoint lines called axes.
%
%
%
%
%

\begin{theorem} \Label{191}
Let $\Sigma_3$ be the intersection of a hwise-fixed-$5$-space and  an  $\HI$-$5$-space, and
denote the set of all fixed-\I-lines  and fixed-\II-lines in $\Sigma_3$ by $\F$.   \begin{enumerate}
\item If $\qtwo$, then  the lines in $\F$ form
  an elliptic  linear congruence.
\item If $q\equiv 0\pmod 3$, then the lines in  $\F \cup \ell$ form a parabolic linear congruence with axis $\ell$ where $\ell$ is the unique ptwise-fixed line in $\Sigma_3$.
  \item If $q\equiv 1\pmod 3$, then the lines in $\F$ form a hyperbolic  linear congruence whose axes are the two ptwise-fixed lines in $\Sigma_3$.
  \end{enumerate}
 \end{theorem}

  \begin{proof} Let $\Sigma_3$ be the intersection of a hwise-fixed-$5$-space and  an  $\HI$-$5$-space, and
denote  the set of fixed-\I-lines  and fixed-\II-lines in $\Sigma_3$ by $\F$. By Theorem~\ref{333}, the points in $\Sigma_3$ are either fixed, $\I_\line$- or $\II_\line$-points. Note that if $P$ is an $X_\line$-point, $X\in\{\I,\II\}$, then $P,P^\sigma,P^{\sigma^2}$ are collinear, and so $P$ lies on a unique fixed-$X$-line. Hence   if two lines from $\F$ meet, then they meet in a fixed point; and
 every non-fixed point of $\Sigma_3$ lies on a unique line of $\F$. Thus
the lines in $\F$ partition the non-fixed points of $\Sigma_3$.

By Lemma~\ref{345b} and Corollary~\ref{221}, $\Sigma_3$ contains exactly $q+1$ fixed-\I-lines which form a 1-regulus $\R$.  The opposite regulus  $\R'$  contains the $\m +1$ ptwise-fixed lines  of $\Sigma_3$.

Suppose $\qzero$, denote the lines of $\R$ by $m_0,\ldots,m_q$ and the unique ptwise-fixed line of $\Sigma_3$ by $\ell$, so $\ell$ is a line of $\R'$.
 Consider the plane $\pi_i=\langle m_i,\ell\rangle$ for any $i\in\{0,\ldots,q\}$. As $\pi_i$ contains two fixed lines, it is fixed. So $\sigma$ induces a projectivity  acting on $\pi_i$ that has $q+1$ fixed points, so has $q+1$ fixed lines, these are necessarily the lines through the point $m_i\cap\ell$. Hence $q-1$ of the lines in $\pi_i$ through the point $m_i\cap\ell$ are fixed-\II-lines. This accounts for $(q+1)(q-1)$ fixed-\II-lines of $\Sigma_3$. This is all the fixed-\II-lines of $\Sigma_3$ by Corollary~\ref{221}.
 That is, for each line $m_i$ in the 1-regulus $\R$, the pencil of lines through $m_i\cap\ell$ in the plane
  $\pi_i=\langle \ell,m_i\rangle$ lies in $\F$, and conversely, each line in $\F$ is one of these lines.
As the planes $\pi_i$ pairwise meet in $\ell$, $\F\cup\ell$ is a set of $q^2+q+1$ lines which form a parabolic congruence with axis $\ell$.

 Suppose $\qone$, so there are
 two  ptwise-fixed lines in $\R'$ which we denote by $\ell_1,\ell_2$. Note that $\Sigma_3$ contains
exactly $2(q+1)$ fixed-points, namely the points on $\ell_1$ and $\ell_2$.
  By Theorem~\ref{294},  each fixed-\II-line in $\Sigma_3$ contains two fixed-points. Hence every line in $\F$ meets both $\ell_1$ and
$\ell_2$. By Corollary~\ref{221}, $\F$ contains $q+1$ fixed-\I-lines (which form a 1-regulus) and $q^2+q$ fixed-\II-lines.
That is, $\F$ is the set of $(q+1)^2$ lines joining a point of $\ell_1$ to a point of $\ell_2$  and so $\F$ is a hyperbolic congruence with axes $\ell_1,\ell_2$.

   Suppose $\qtwo$, then by Corollary~\ref{221},  $\Sigma_3$ contains no fixed
points, so the $q^2+1$ lines in $\F$ partition the points of $\Sigma_3$, that is,
$\F$ is a 1-spread. We show $\F$ is regular by looking in the quadratic extension.
 Embed $\PG(2,q^3)$ in $\PG(2,q^6)$ and consider the Bose representation of $\PG(2,q^6)$ in $\PG(8,q^2)$.
 As $\qtwo$, in the quadratic extension we have $q^2 \equiv 1 \pmod{3}$. Hence by the preceding paragraph,  in the quadratic  extension of $\Sigma_3$ to a 3-space over $\Fqq$,   the fixed-\I-lines and fixed-\II-lines form a hyperbolic congruence with axes two   ptwise-fixed lines denote $\ell,\ell'$. These axes lie in $\PG(8,q^2)\setminus\PG(8,q)$ (ie are `imaginary').
The hyperbolic congruence must contain all the fixed-\I-lines and fixed-\II-lines in $\F$ from $\PG(8,q)$ (ie the `real'  fixed lines). That is, the $\Fqq$-extension of the real fixed-\I-lines and fixed-\II-lines in $\F$ each meet both $\ell$ and $\ell'$. Hence $\ell$ and $\ell'$ are transversals of the spread $\F$. That is, the  spread $\F$   is a regular spread, and so is an elliptic congruence.
   \end{proof}

    We can extend Theorem~\ref{191} to describe the   fixed-\I-lines and fixed-\II-lines in a hwise-fixed-$5$-space. First we look at the
set of fixed-\I-lines.

   \begin{corollary}\Label{333d}
The   $q^2+q+1$ fixed-\I-lines in a  hwise-fixed-$5$-space are the ruling
lines of a Segre variety $\mathcal S_{1;2}$. The $q+1$ ruling planes of this variety are all $\piq$-ruling-planes, and include $\m+1$ ptwise-fixed planes.
\end{corollary}

\begin{proof} Let $\Pifix$ be a  hwise-fixed-$5$-space. By Theorem~\ref{333}, $\Pifix$ contains $q^2+q+1$ fixed-\I-lines. By
Lemma~\ref{900}, each lies in a distinct $\SI$-plane. As discussed in Section~\ref{sec2.2}, the
set of all $\SI$-planes is  $\Bo{\piq}$, and forms one set of ruling planes of a Segre variety $\S_{2;2}$. The other set of ruling planes are called  $\piq$-ruling-planes, and include the $\g$ ptwise-fixed planes by Theorem~\ref{000}. So the set of fixed-\I-lines in $\Pifix$ is the intersection of a $\S_{2;2}$ with a $5$-space, so is a Segre variety $\S_{1;2}$. By Lemma~\ref{294}, a fixed-\I-line contains $\m+1$ fixed points, so the  $\S_{1;2}$ Segre variety contains $\m+1$ ptwise-fixed planes.
\end{proof}

{\bfseries Remark.\ } Let  $\L$  be the set of all fixed-\I-lines and fixed-\II-lines in a hwise-fixed-$5$-space $\Pifix$, then the
 lines of $\L$ partition the non-fixed-points of $\Pifix$. Theorem~\ref{191} can be generalised to show the following.
 \begin{itemize}
  \item If $\qtwo$, then $\L$ is a regular 1-spread of the $5$-space $\Pifix$.
\item
 If $\qzero$, then  the $5$-space $\Pifix$ contains one ptwise-fixed plane $\pifix$, and
each line of $\L$ meets $\pifix$;
moreover each $3$-space of $\Pifix$ containing $\pifix$ contains $q^2$ lines of $\L$ through a common vertex.
 \item  If $\qone$, then the $5$-space $\Pifix$ contains two ptwise-fixed planes, and $\L$ is the set of all lines which meet both ptwise-fixed planes.
 \end{itemize}

 \subsection{Planes fixed by $\sigma$}

 In this section, we describe the planes of $\PG(8,q)$ which are fixed by
$\sigma$. We know that $\SI$-planes are fixed. In Theorem~\ref{000}, we showed that there are exactly $\g$ ptwise-fixed planes, moreover these are $\piq$-ruling-planes. We now show that there are three further types of fixed planes.

 \begin{theorem}\Label{112}
  Let $\pi$ be a fixed plane of $\PG(8,q)$, then $\pi$ is exactly one of the following types.
\begin{itemize}
\item  $\pi$  is  a ptwise-fixed plane;
\item $\pi$ is an $\SI$-plane;

\item   $\pi$   contains at least one $\II_\plane$-point, in which case $\pi$ lies in a  unique $\HI$-$5$-space, and we call this   a \emph{fixed-\II-plane;}
\item  $\pi$   contains at least one $\III_\plane$-point,  we call this
a  \emph{fixed-\III-plane; }
\item $\pi$ is contained in a hwise-fixed-$5$-space and $\pi$ is not fixed pointwise, we call this   an  \emph{\h-plane}.
 \end{itemize}
 \end{theorem}

\begin{proof}
Let $\pi$ be a  plane of $\PG(8,q)$ which is fixed by $\sigma$, but is not fixed pointwise, and  $\pi$ is  \emph{not} contained in a hwise-fixed-$5$-space.
As $\pi$ is not fixed pointwise, $\sigma$ induces a non-identity projectivity  acting on $\pi$, so there are   at most $q+2$ fixed points in $\pi$.
Suppose  the non-fixed points of $\pi$ are all either $\I_\line$ or $\II_\line$-points, then by
 Theorem~\ref{333a}, every point of  $\pi$ is contained in a hwise-fixed-$5$-space, which is a contradiction.
Hence
$\pi$ contains a point $P$ which is a $Y_\plane$-point for some $Y\in\{\I,\II,\III\}$, so $\pi=\langle P,P^\sigma,P^{\sigma^2}\rangle$.

If $P$ is  a $\I_\plane$-point, then  $\pi$ is an $\SI$-plane.
 If $P$ is  a  $\II_\plane$-point,  then $P,P^\sigma,P^{\sigma^2}$ lie in
a common $\HI$-$5$-space, so $\pi$ is contained in  an  $\HI$-$5$-space and so $\pi$ contains only \I-points and \II-points. We call $\pi$ a fixed-\II-plane, and note that it does not contain a \III-point. The final possibility for $\pi$ is if $P$ is a
 $\III_\plane$-point, in which case the above argument shows that $\pi$ does not contain a $\II_\plane$-point. Hence we have the five disjoint classes given in the result.
 \end{proof}

So there are five types of fixed planes in $\PG(8,q)$. We have previously classified the  points and fixed lines  in an $\SI$-plane and in a ptwise-fixed plane. We next classify the points and fixed lines in the other three classes of fixed planes.

 \begin{lemma}\Label{555}
Let $\pi$ be a  fixed-$Y$-plane for some $Y\in\{\II,\III\}$. Then $\pi$ contains exactly $\g$ fixed points and    $q^2+q+1-\g (q-\m +1)$ $Y_\plane$-points; and exactly $\g$ fixed lines .
\end{lemma}

 \begin{proof}  Let $\pi$ be a  fixed-$Y$-plane for some $Y\in\{\II,\III\}$, so
 $\pi=\langle P,P^\sigma,P^{\sigma^2}\rangle$ for some $Y_\plane$-point
$P$.  Suppose $\pi$ contains a line $m$ which is fixed pointwise by $\sigma$, then $\sigma$ induces a collineation acting on $\pi$ that fixes  a point $Q$ linewise.   Hence $Q,P,P^\sigma,P^{\sigma^2}$ are collinear, a contradiction. Hence $\pi$ meets a ptwise-fixed plane  in at most one point.

 If $\g =1$, then $q^2+q+1\equiv 1\pmod3$, so $\sigma$ fixes at least one point of $\pi$. Thus  $\pi$ contains exactly one fixed point.

If  $\g =3$, then $q^2+q+1\equiv0 \pmod3$. By Theorem~\ref{000},
there are three ptwise-fixed planes, so $\pi$ contains either zero or three fixed points (forming a triangle).
Let $\X$ be the set of planes of the latter type, that is, 
the planes in $\X$ contain one point from each of the three ptwise-fixed planes. So $|\X|=(q^2+q+1)^3$. Let $\beta$ be a plane in $\X$, then
$\beta$ contains exactly three fixed points, hence exactly  three fixed lines, and so $q^2-2q+1$ $Z_\plane$-points for some $Z\in\{\I,\II,\III\}$.
As the planes in $\X$ are fixed, a $Z_\plane$-point lies in at most one plane of $\X$.
Hence the planes in $\X$ cover exactly $k=(q^2+q+1)^3(q^2-2q+1)$ $Z_\plane$-points.
We claim that $k$ is equal to the total number of $Z_\plane$-points. Then $\pi$ cannot contain zero fixed points. To verify the claim, note that the number of $\SI$-planes, $\SII$-planes and $\SIII$-planes is $s_1=q^2+q+1$, $s_2=(q^2+q+1)(q^3-q)$, and  $s_3=q^3(q-1)^2(q+1)$  respectively. By  Corollary~\ref{294b}, the number of  $\I_\plane$-points in an $\SI$-plane  is  $a_1=q^2-2q+1$.
By  Corollary~\ref{333b1}, the number of  $\II_\plane$-points in an $\SII$-plane  is   $a_2=q^2+q-2$; and
 the number of  $\III_\plane$-points in an $\SIII$-plane  is   $a_3=q^2+q+1$.   It is straightforward to verify that $k=s_1a_1+s_2a_2+s_3a_3$. Hence every $Z_\plane$-point, $Z\in\{\I,\II,\III\}$,  lies in exactly one plane of $\X$ as required.
 We conclude that $\pi$ contains exactly three fixed points (which are non-collinear) one from each ptwise-fixed plane.

We have shown that $\pi$ contains exactly  $\g $  fixed points. As $\pi$ is fixed, $\sigma$ induces a collineation on $\pi$ which fixes exactly  $\g $ lines. These lines are not fixed ptwise as $\pi$ only contains $\g$ fixed points.
So the $\g$ fixed lines in $\pi$ are  fixed-$X$-lines  for some $X\in\{\I,\II\}$.
Hence by Theorem~\ref{294}, $\pi$ contains: $\g $ fixed points; $\g(q-\m)
$ points that are $\I_\line$- or $\II_\line$-points. By Theorem~\ref{112}, the remaining $q^2+q+1-\g (q-\m +1)$ points are all either  $\II_\plane$-points or  $\III_\plane$-points.  \end{proof}

 \subsubsection{Fixed-\II-planes}

\begin{theorem}\Label{dark-b}
There are two types of fixed-\II-planes.
 \begin{enumerate}
 \item A \emph{fixed-\IIa-plane} contains
\begin{itemize}
\item  $\g $ fixed points, $q-\m $ $\I_\line$-points,
  $(\g -1)(q-\m)$ $\II_\line$-points and  $q^2+q+1-\g (q-\m +1)$ $\II_\plane$-points;
\item  $1$ fixed-\I-line and $\g -1$ fixed-\II-lines.
\end{itemize}
The total number of fixed-\IIa-planes is $\g (q-2+\g-\m)(q+1)(q^2+q+1)$.
 \item A \emph{fixed-\IIb-plane} contains
 \begin{itemize}
 \item  $\g $ fixed points, $\g (q-\m )$ $\II_\line$-points and $q^2+q+1-\g (q-\m +1)$ $\II_\plane$-points;
 \item   $\g $ fixed-\II-lines
 \end{itemize}
The total number of fixed-\IIb-planes is $(q^3-q)(q^2+q+1)$.
\end{enumerate}
\end{theorem}

\begin{proof} Let $\pi$ be a fixed-\II-plane,  so by Lemma~\ref{555}, $\pi$ contains exactly  $\g $  fixed points and $\g$  fixed-$X$-lines  for some $X\in\{\I,\II\}$.
By Lemma~\ref{345}, a fixed-\I-line lies in an $\SI$-plane. As $\pi$ is not an $\SI$-plane, it contains  at most  one fixed-\I-line, giving the two possibilities described in the result.

We next count the number of planes of each type. By Theorem~\ref{112}, a fixed-\II-plane is contained in  a  unique $\HI$-$5$-space, so we count the fixed-\II-planes contained in one $\HI$-$5$-space and multiply by $q^2+q+1$ (the number of $\HI$-$5$-spaces).

Suppose $\qtwo$, and let $\Sigma_5$ be an $\HI$-$5$-space.
By Corollary~\ref{333e}, $\Sigma_5$ contains exactly one ptwise-fixed line (denoted $\ell$), $q+1$
fixed-\I-lines and $q^2-q$ fixed-\II-lines. By Corollary~\ref{294b}, $\Sigma_5$ contains exactly $q+1$ fixed points, so all the fixed points in $\Sigma_5$ lie on $\ell$.
By Theorems~\ref{333a} and~\ref{191}, the set of fixed-\I-lines and fixed-\II-lines of $\Sigma_5$ form a spread in the $3$-space
$\Pi_3=\Pifix\cap\Sigma_5$. Moreover by  Theorem~\ref{294},  $\ell$ does not meet any fixed-$X$-line, $X\in\{\I,\II\}$. Hence the fixed lines in $\Sigma_5$ are pairwise disjoint.
Let $\pi$ be a fixed plane contained in $\Sigma_5$.
As $q^2+q+1 \equiv 1 \pmod{3}$ and $\sigma$ has order 3, $\pi$ contains at least one fixed point and hence contains at least one fixed line. As the fixed lines in $\Sigma_5$ are pairwise disjoint, $\pi$ contains exactly one fixed point $P\in\ell$,  and exactly one fixed line $m\subset\Pi_3$.
There are three possibilities for $\pi$:
(a) if $m$
 is a fixed-\I-line and $P,m$ lie in the same $\SI$-plane, then $\pi$ is a $\SI$-plane;
(b) if $m$ is a fixed-\I-line and $P,m$ lie in different $\SI$-planes, then $\pi$ is a fixed-\IIa-plane;
(c) if $m$ is a fixed-\II-line, then $\pi$ is a fixed-\IIb-plane.
So the numbers of such planes in $\Sigma_5$ are
$(q+1)$, $(q+1)q$ and $(q^2-q)(q+1)$
respectively. Multiplying by the number of $\HI$-$5$-spaces gives the  number of fixed-\II-planes as stated in the lemma.

Suppose $\qzero$, and let $\Sigma_5$ be an $\HI$-$5$-space.
By Corollary~\ref{333e}, $\Sigma_5$ contains exactly one ptwise-fixed line denoted $\ell$; $q+1$ fixed-\I-lines denoted $m_0,\ldots,m_q$ which form a 1-regulus; and $q^2-1$ fixed-\II-lines. By Corollary~\ref{294b}, $\Sigma_5$ contains exactly $q+1$ fixed points, so all the fixed points in $\Sigma_5$ lie on $\ell$.
By Theorems~\ref{333a} and~\ref{191},  the fixed-\I-lines and fixed-\II-lines of $\Sigma_5$ form a parabolic congruence with axis $\ell$ in the 3-space $\Pi_3=\Pifix\cap\Sigma_5$. That is, the
fixed lines of $\Sigma_5$ are  the set of lines  that contain   the point
$m_i\cap\ell$ and lie in the plane
 $\pi_i=\langle m_i,\ell\rangle$, for some $i\in\{0,\ldots,q\}$.

 Let $\pi$ be a fixed plane in $\Sigma_5$. Consider the possibilities for the fixed subspace $\pi\cap\Pi_3$.
 If $\pi$ lies in $\Pi_3$, then $\pi$ is contained in a hwise-fixed-$5$-space and so is  an  $\h$-plane. If $\pi\cap\Pi_3$ is a point, then it is a fixed point, hence $\pi$ contains a fixed line, contradicting all fixed lines of $\Sigma_5$ lying in $\Pi_3$. If $\pi\cap\Pi_3$ is the line $\ell$, then $\pi$ contains $q+1$ fixed points, hence contains $q+1$ fixed lines,  contradicting all fixed lines of $\Sigma_5$ lying in $\Pi_3$.  Hence a fixed plane is either contained in $\Pi_3$ or meets $\Pi_3$ in a fixed-$X$-line, $X\in\{\I,\II\}$.

 Let $t$ be a fixed-$X$-line, $X\in\{\I,\II\}$, so $t$ lies in $\pi_i$ for some $i\in\{0,\ldots,q\}$. Suppose $t$ lies in a fixed plane $\pi$ with $\pi$ not contained in $\Pi_3$. Then the fixed 3-space $\Sigma_3=\langle \pi, \pi_i\rangle$ contains exactly $q+1$ fixed points on the line $\ell$, hence contains $q+1$ fixed planes through the  line $t$. If $t$ lies on a further fixed plane $\pi'$, with $\pi'\not\subset\Sigma_3$,  $\pi'\not\subset\Pi_3$, then the fixed 3-space $\langle \pi,\pi'\rangle$ contains exactly one fixed point, and at least two fixed planes, a contradiction. Hence $t$ lies in at most $q+1$ fixed planes. As $\Pi_3$ is fixed and contains exactly $q+1$ fixed points, it contains exactly $q+1$ fixed planes, namely $\pi_0,\ldots,\pi_q$. Hence the number of fixed planes containing a fixed-\I-line or a fixed-\II-line is at most $x=(q^2+q)q+(q+1)$.

We now show that this bound is attained by counting the total number of fixed planes in $\Sigma_5$.
Consider in turn the five different types of fixed planes described in Theorem~\ref{112}. Firstly, by Theorem~\ref{000}, the ptwise-fixed-plane is a $\piq$-ruling plane and is not contained in $\Sigma_5$. Secondly, there are exactly  $q+1$ $\SI$-planes in $\Sigma_5$. Thirdly, each $\II_\plane$-point in $\Sigma_5$ lies in a unique fixed-\II-plane. By Corollary~\ref{333b1}, the number of $\II_\plane$-points in $\Sigma_5$ is  $a=(q^3-q)(q^2+q)$. By Theorem~\ref{dark-b}, each fixed-\II-plane contains exactly $b=q^2$ $\II_\plane$-points. Hence there are $a/b=q^3+q^2-q-1$ fixed-\II-planes in $\Sigma_5$. Next, there are no $\III$-points in $\Sigma_5$,  so there are no fixed-\III-planes in $\Sigma_5$. Finally, as argued above, there are $q+1$ \h-planes in $\Sigma_5$, namely $\pi_0,\ldots,\pi_q$. Hence $\Sigma_5$ contains exactly $y=q^3+q^2+q+1$ fixed planes. Moreover, each of these planes contains at least one fixed-$X$-line for some $X\in\{\I,\II\}$. As  $x=y$, each fixed-$X$-line in $\Pi_3$ lies in exactly $q+1$ fixed planes of $\Sigma_5$.
Hence each fixed-\I-line in $\Pi_3$ lies in one \h-plane, one $\SI$-plane and $q-1$ fixed-\IIa-planes of $\Sigma_5$, so $\Sigma_5$ contains exactly $(q-1)(q+1)$ fixed-\IIa-planes.
Further, each fixed-\II-line in $\Pi_3$ lies in one \h-plane and $q$ fixed-\IIb-planes of $\Sigma_5$, so $\Sigma_5$ contains exactly $q(q^2-1)$ fixed-\IIb-planes. Multiplying by the number of $\HI$-$5$-spaces gives the  number of fixed-\II-planes as stated in the lemma.

Suppose $\qone$, and let $\Sigma_5$ be an $\HI$-$5$-space.
By Corollary~\ref{333e}, $\Sigma_5$ contains exactly three ptwise-fixed lines, denote these $\ell_1,\ell_2,\ell_3$. By Lemma~\ref{555}, each  fixed-\II-plane contains three fixed points, one from each of $\ell_1,\ell_2,\ell_3$.
Let $\X$ be the set of  $(q+1)^3$ planes which meet $\ell_1$, $\ell_2$ and $\ell_3$. If $\pi$ is a plane in $\X$, then
$\pi$ contains one point from each ptwise-fixed-plane, so by Theorem~\ref{333}, $\pi$ is not contained in a hwise-fixed-$5$-space. As $\Sigma_5$ contains no $\III$-points, by Theorem~\ref{112},
 $\pi$ is either an $\SI$-plane or a fixed-\II-plane.
 Let $\alpha$ be an $\SI$-plane contained in $\Sigma_5$. By Corollary~\ref{294b}, $\alpha$ contains three fixed-\I-lines, denote these $m_1,m_2,m_3$ such that $m_i\cap m_j\in\ell_k$ with $\{i,j,k\}=\{1,2,3\}$. The fixed-\I-line $m_1$ lies in exactly $q$ fixed-\IIa-planes, namely the planes $\langle m_1,X\rangle$ for $X\in\ell_1$, $X\notin\alpha$. Further, each of these fixed-\IIa-planes contains exactly one fixed-\I-line, namely $m_1$. There are $q+1$ choices for the $\SI$-plane $\alpha$ and $3$ choices for the fixed-\I-line $m_1$ in $\alpha$. Hence the number of fixed-\IIa-planes in $\Sigma_5$ is $3q(q+1)$. Thus the number of fixed-\IIb-planes in $\Sigma_5$ is $|\X|-(q+1)-3q(q+1)=q^3-q$.
 Multiplying by the number of $\HI$-$5$-spaces gives the  number of fixed-\II-planes as stated in the lemma.
\end{proof}

 \subsubsection{Fixed-\III-planes}\Label{764}

\begin{theorem}\Label{765}
The number of fixed-\III-planes is $q^3(q^2-1)(q-1)(q^2+q+1)/(q^2+q+1-\g (q-\m +1))$.
Each fixed-\III-plane contains exactly
\begin{itemize}
\item $\g $ fixed points, $\g (q-\m )$ $\II_\line$-points and
$q^2+q+1-\g (q-\m +1)$
$\III_\plane$-points;

\item
$\g $ fixed-\II-lines.
 \end{itemize}
\end{theorem}

\begin{proof}
Let $\pi$ be a fixed-\III-plane,  so by  Lemma~\ref{555}, $\pi$ contains exactly  $\g $  fixed points and $\g$  fixed  lines.
 Suppose $\pi$ contains a
 fixed-\I-line $\ell$.
 Let $\alpha$ be the $\SI$-plane  containing $\ell$.
Let $C$ be a $\III$-point in $\pi$ and let $\gamma$ be the unique $\S$-plane  containing $C$. Note that $C$ is a $\III_\plane$-point  and $\gamma$ is an $\SIII$-plane.
The $5$-space $\Sigma_5=\langle \alpha,\gamma\rangle$ contains $\pi$. As $\alpha,\gamma$ are $\S$-planes, $\Sigma_5$ is
 an $\H$-$5$-space.
 As $\pi$ is fixed, $C^\sigma,C^{\sigma^2}$ lie in $\pi$, so they lie in $\Sigma_5$. As
 $\Sigma_5$ contains $C^\sigma,C^{\sigma^2}$ and $\Sigma_5$ is an $\H$-$5$-space, $\Sigma_5$ also
  contains the $\S$-planes containing $C^\sigma,C^{\sigma^2}$. That is, $\Sigma_5$ contains $\gamma,\gamma^\sigma,\gamma^{\sigma^2}$, which contradicts the definition of a Type-III-point in $\PG(2,q^3)$. Hence $\pi$ contains $\g $ fixed-\II-lines and so by Theorem~\ref{294}, $\pi$ contains  $\g (q-\m )$ $\II_\line$-points. The remaining
$q^2+q+1-\g (q-\m +1)$ points are
$\III_\plane$-points.

 To count the fixed-\III-planes, recall that there
are $ q^3(q^2-1)(q-1)$ $\SIII$-planes, and so  there are $ a=q^3(q^2-1)(q-1)(q^2+q+1)$ $\III_\plane$-points in $\PG(8,q)$. By Lemma~\ref{765}, the number of $\III_\plane$-points in a fixed-\III-plane is  $b=q^2+q+1-\g (q-\m +1)$. As each $\III_\plane$-point lies in exactly one fixed plane,
the number of fixed-\III-planes in $\PG(8,q)$ is $a/b$ as required.
   \end{proof}

We conclude this section with a result on the structure of fixed-\III-planes.

\begin{theorem}\Label{401}
Let $\gamma$ be an $\SIII$-plane. The set of all fixed-\III-planes that meet $\gamma$ form one system of ruling planes of a Segre variety $\S_{2;2}$.
\end{theorem}

\begin{proof} Let $\gamma$ be an $\SIII$-plane. By Corollary~\ref{333b1}, each point in $\gamma$ is a  $\III_\plane$-point. Hence the set of fixed-\III-planes meeting $\gamma$ is    $\mathbb G=\{\langle G,G^\sigma,G^{\sigma^2}\rangle\,|\,G\in\gamma\}$ and has size $q^2+q+1$.
Moreover, each plane in $\mathbb G$ meets $\gamma^\sigma$ in a point and meets $\gamma^{\sigma^2}$ in a point. As the three planes $\gamma,\gamma^\sigma,\gamma^{\sigma^2}$ span $\PG(8,q)$, it follows from \cite[Corollary 2.2]{BJW3} that the planes in $\mathbb G$ are pairwise disjoint.
By Theorem~\ref{000}, there are $\g$ ptwise-fixed planes, denote one of these by $\pifix$. It follows from Theorem~\ref{765} that each plane in $\mathbb G$ meets $\pifix$ in a distinct point.

Consider a plane $\pi\in\mathbb G$, so $\pi=\langle G,G^\sigma,G^{\sigma^2}\rangle$ for some point $G\in\gamma$, and $\pi$ contains a $\I$-point of $\pifix$, denoted $A$. As $A$ is fixed, it does not lie on a side of the triangle $GG^\sigma G^{\sigma^2}$. It follows that the four planes $\pifix,\gamma,\gamma^\sigma,\gamma^{\sigma^2}$ are such that any three span $\PG(8,q)$. Hence they are ruling planes of a   unique Segre variety $\S_{2;2}$, and the second
  system of ruling planes of this Segre variety $\S_{2;2}$ is $\mathbb G$.
  \end{proof}

 \subsubsection{Fixed \h-planes}

We now look at the final class of fixed planes given in Theorem~\ref{112}, namely \h-planes which are fixed planes that are contained in a hwise-fixed-$5$-space and are not ptwise-fixed. We show that these do not occur if $\qtwo$, otherwise there are two types of \h-planes.

\begin{theorem}\Label{113}
\begin{enumerate}
\item If $\qtwo$, then there are no \h-planes.

\item\Label{b1} If $q\not\equiv-1\pmod3$, then there are two types of \h-planes.
\begin{enumerate}
\item \Label{b2a}
An \ha-plane lies in an $\HI$-$5$-space and contains
\begin{itemize}
\item
 $q+1+\m$ fixed-points,    $q-\m$ $\I_\line$-points  and
 $q^2-q$ $\II_\line$-points,

 \item
 $1$ ptwise-fixed line, $1$ fixed-\I-line and
 $q-1+\m$ fixed-\II-lines.
 \end{itemize}
There are
$\g(\m+1)(q+1)(q^2+q+1)$ $\ha$-planes.
 \item\Label{b2b} An \hb-plane contains
 \begin{itemize}
 \item
 $q+1+\m$ fixed-points  and
 $q^2-\m$ $\II_\line$-points,

 \item
$1$ ptwise-fixed line and  $q+\m$ fixed-\II-lines.

  \end{itemize}
There
   are  $\g(\m+1)(q^2+\m-1)(q^2+q+1)$ \hb-planes.
   \end{enumerate}
\end{enumerate}
\end{theorem}

\begin{proof} Let $\Pifix$ be a hwise-fixed-$5$-space and suppose $\pi$ is an \h-plane   contained in $\Pifix$. By Theorem~\ref{333}, the points of $\pi$ are either  fixed points,  $\I_\line$-points or  $\II_\line$-points.
As $\pi$ is fixed, $\sigma$ induces a non-identity projectivity  on $\pi$, so the number of fixed points of $\pi$ is either $0$, $1$, $2$, $3$, $q+1$
or $q+2$.  Let $P$ be a non-fixed point in $\pi$, so $P$ is an $X_\line$-point with $X\in\{\I,\II\}$, and so $\langle P,P^\sigma,P^{\sigma^2}\rangle$ is a  fixed-$X$-line contained in $\pi$.   Suppose $\pi$ contains two fixed-\I-lines $m_1,m_2$, then $m_1,m_2$ meet, so lie in the same $\SI$-plane  by Theorem \ref{294}. Hence $\pi$ is an $\SI$-plane contained in $\Pifix$, contradicting Lemma~\ref{900}. Hence $\pi$ contains at most one fixed-\I-line.
 So by Theorem~\ref{294}, $\pi$ contains at most $q-\m $ $\I_\line$-points.
Hence $\pi$ contains at least $q^2+q+1-(q+2)-(q-\m )
=q^2-q-1+\m $ $\II_\line$-points. Each of these $\II_\line$-points lies
on a unique fixed-\II-line which is contained in $\pi$ (as $\pi$ is fixed). By Theorem~\ref{294}, a fixed-\II-line contains $q-\m $ $\II_\line$-points, so $\pi$ contains at least $q-2$ fixed-\II-lines. Dually, $\pi$ contains at least $q-2$
fixed points. Hence $\pi$ contains either $q+1$ or $q+2$ fixed points.
Note that this means that $\Pifix$ contains at least $q+1$ fixed points, which contradicts Theorem~\ref{333} if $\qtwo$.
Hence there are no \h-planes when $\qtwo$.

As $\pi$ contains $q+1$ or $q+2$ fixed points, it contains $q+1$ or $q+2$
fixed lines, one of which is a ptwise-fixed line, denote this by $\ell$.
Suppose $\pi$ contains a fixed-\I-line $m$, so $\pi=\langle\ell,m\rangle$.  By Theorem~\ref{000}, $\ell$ meets $q+1$ $\SI$-planes, these lie in a unique $\HI$-$5$-space denoted $\Sigma_5$.
Hence the fixed line $m$ lies in one of the $\SI$-planes that $\ell$ meets. That is, $m$
lies in $\Sigma_5$, and so $\pi$ lies in $\Sigma_5$. That is, an \h-plane
which contains a fixed-\I-line lies in  an  $\HI$-$5$-space.

Suppose $\qzero$. As $q^2+q+1\equiv 1\pmod3$ and $\sigma$ has order 3, the number of fixed points in $\pi$ is congruent to 1 modulo 3. Hence $\sigma$ fixes exactly $q+1$ points of $\pi$. That is, $\sigma$
 induces an elation on $\pi$ with axis the ptwise-fixed line $\ell$, and centre a point $P\in\ell$. Recall that $\pi$ contains at most one fixed-\I-line, so the lines of $\pi$ through $P$ consist of $\ell$, at most one fixed-\I-line and at least $q-1$ fixed \II-lines.
 The number of points of each type in an \h-plane now follows from Theorem~\ref{294}.

If $\qone$, then a similar argument shows that $\sigma$ induces a homology on $\pi$ with axis $\ell$ and centre $P\notin\ell$.
 By Theorem~\ref{294}, each fixed-$X$-line, $X\in\{\I,\II\}$,  contains exactly two fixed points, so the points and fixed lines in $\pi$ are as described.

 This completes the description of the two types of \h-planes, we now count them.
Suppose $\qzero$ and let $\Pifix$ be the unique hwise-fixed-$5$-space. By Theorem~\ref{333}, $\Pifix$ contains the unique ptwise-fixed plane $\pifix$, and $q^2+q+1$ fixed-\I-lines and $q^4+q^3-q-1$ fixed-\II-lines. We proved above that each \h-plane meets $\pifix$ in a line, and contains $q$ fixed-$X$-lines where $X\in\{\I,\II\}$.
We count the $\ha$-planes containing a fixed-\I-line $m$ as follows. The line $m$  meets $\pifix$ in a point, so meets $q+1$ lines of $\pifix$, denoted $\ell_i$, $i=0,\ldots,q$. This gives $q+1$ $\h$-planes $\langle m,\ell_i\rangle$ containing $m$. Further, every \h-plane containing $m$ meets $\pifix$ in a line. Hence there are $(q+1)(q^2+q+1)$ \ha-planes. We use a similar technique to count $x$, the total number of
 \h-planes in $\Pifix$. We count the ordered pairs $(m,\pi)$ where $m$ is a fixed-$X$-line, $X\in\{\I,\II\}$, and $\pi$ is an \h-plane containing $m$. We have $(q^4+q^3+q^2)(q+1)=qx$, hence $x=(q^2+q)(q^2+q+1)$. As $(q+1)(q^2+q+1)$ of these are \ha-planes, the number of \hb-planes is $(q^2-q)(q^2+q+1)$, as required.

Suppose $\qone$ and let $\Pifix$ be one of the three hwise-fixed-$5$-spaces. By Theorem~\ref{333}, there are two ptwise-fixed planes $\pifix,\pifix'$ contained in $\Pifix$.  We proved above that  the \h-planes in $\Pifix$ are precisely the planes that meet one of $\pifix,\pifix'$ in a  point and meet the other in a line. Hence there are $2(q^2+q+1)^2$ \h-planes in $\Pifix$. By Corollary~\ref{333d}, there are $q^2+q+1$ fixed-\I-lines in $\Pifix$, and they are pairwise disjoint. We count the $\ha$-planes containing a fixed-\I-line $m$ as follows. The line $m$  meets $\pifix$ in a point, so meets $q+1$ lines of $\pifix$, denoted $\ell_i$, $i=0,\ldots,q$. This gives $q+1$ $\h$-planes $\langle m,\ell_i\rangle$ containing $m$. By symmetry, we have $q+1$ \ha-planes containing $m$ and meeting $\pifix'$ in a line. Hence $\Pifix$ contains $2(q+1)(q^2+q+1)$ $\ha$-planes and so contains $2(q^2+q+1)^2-2(q+1)(q^2+q+1)$ \hb-planes. The result follows as there are $\g=3$ hwise-fixed-$5$-spaces.

Note that if $\qtwo$, by part 1, there are no \h-planes. As $\m+1=0$, both the counts in part 2 are  zero in this case, so these counts also hold when $\qtwo$.
\end{proof}

\subsection{Summary tables of counting}\Label{009}

We present three tables to summarise the counting in $\PG(8,q)$ (these appear at the end of the article).
In Table~\ref{table16}, we count the various subspaces. The number of $\SX$-planes, $X\in\{\I,\II,\III\}$,  is straightforward to compute by counting points in $\PG(2,q^3)$. The number of points of each type follows  from the number of $\SX$-planes and Results~\ref{000}, \ref{294b}, \ref{333b1}.
The number of fixed lines of each type follows from Lemmas~\ref{345} and~\ref{333c}. The number of fixed planes of each type follows from Theorems~\ref{dark-b}, \ref{765} and \ref{113}.
In Tables~\ref{table14}  and \ref{table18}  we list respectively the number of points/fixed lines of each type in our named subspaces.  The number
of points of each type in a fixed line  follows from Theorem~\ref{294}. The number of points/fixed lines of each type in an $\S$-plane follows from Corollaries~\ref{294b} and \ref{333b1}. The number of points and fixed lines of each type in a fixed  plane follows from Results~\ref{dark-b}, \ref{765} and \ref{113}.
The number of points/lines of each type in a hwise-fixed subspaces is given in Theorem~\ref{333}.


\section{Related $\phi$-fixed linear sets  of $\PG(2,q^3)$}\Label{Sec3}

Let $\Pi_r$ be an $r$-space of $\PG(8,q)$. We define $\B(\Pi_r)$ to be the set of points of $\PG(2,q^3)$ that correspond to an $\S$-plane which meets $\Pi_r$ in a non-empty subspace. That is,
$$\B(\Pi_r)=\{\bar P\in\PG(2,q^3)\,|\, \Pi_r\cap\Bo{P}\neq \emptyset\}.$$
The set $\B(\Pi_r)$ is called a \emph{linear set of $\PG(2,q^3)$ of rank $r+1$.}
For more information on linear sets, including the definition of clubs and scattered linear sets of $\PG(1,q^3)$, see \cite{fieldreduction}.

A fixed subspace of $\PG(8,q)$ gives rise to a linear set of $\PG(2,q^3)$ that is fixed by $\phi$.
It is straightforward to use the preceding results of this article to describe the  $\phi $-fixed linear sets of $\PG(2,q^3)$ arising from $\sigma$-fixed subspaces of $\PG(8,q)$.

\begin{table}[h!]
\caption{Theorem~\ref{thm:linear}}
\begin{center}
\begin{tabular}{lp{13cm}}
\toprule
$\Pi$& $\B(\Pi)$, the  $\phi $-fixed linear set  of $\PG(2,q^3)$ corresponding to $\Pi$\\
\toprule
fixed point&  a Type-\I-point\\
\midrule
ptwise-fixed line&   an $\Fq$-line  containing $q+1$ Type-\I-points\\
fixed-\I-line& a Type-\I-point\\
fixed-\II-line& an $\Fq$-line contained in a Type-\I-line; contains $\m +1$ Type-\I-points and $q-\m $ Type-\II-points\\
\midrule

ptwise-fixed plane& the $\Fq$-plane $\piq$ (which contains the $q^2+q+1$ Type-\I-points)\\
fixed-\IIa-plane& a club contained in a  Type-\I-line; contains  $\g-\m$ Type-\I-points, one if which is the head, and $q^2-\g+\m+1$ Type-\II-points\\
fixed-\IIb-plane&a scattered linear set contained in a  Type-\I-line; contains $\g$ Type-\I-points and $q^2+q+1-\g$ Type-\II-points\\
fixed-\III-plane& an $\Fq$-plane; contains $\g $ Type-\I-points, $\g (q-\m )$ Type-\II-points and
$q^2+q+1-\g (q-\m +1)$ Type-\III-points\\
\ha-plane& a club contained in a  Type-\I-line, with $q+1$ Type-\I-points (one of which is the head) and $q^2-q$ Type-\II-points; \\
\hb-plane &either: a scattered linear set contained in a Type-\I-line (if the \hb-plane is in an $\HI$-$5$-space); or otherwise, an $\Fq$-plane, with $q+1+\m$ Type-\I-points (on an $\Fq$-line, see Thm \ref{113}) and  $q^2-\m$ Type-\II-points\\

\midrule
hwise-fixed-$5$-space& a linear set of $\PG(2,q^3)$ of rank 6 which is the set of all Type-\I- and Type-\II-points of $\PG(2,q^3)$.\\
\bottomrule
\end{tabular}
\end{center}
\label{tableB}
\end{table}

\begin{theorem} \label{thm:linear}
Let $\Pi$ be a subspace of $\PG(8,q)$ that is fixed by $\sigma$. Then $\B(\Pi)$ is a linear set of $\PG(2,q^3)$ that is fixed by $\phi$. The description of $\B(\Pi)$ for the different cases when $\Pi$ is a fixed point,  fixed line,  fixed plane or a hwise-fixed-$5$-space is given in
 Table~\ref{tableB}.
\end{theorem}

One interesting distinction related to the value of $\g=\gcd(3,q-1)$ is that if  $\g=1$, then  the  only  $\phi$-fixed-$\Fq$-subplanes of $\PG(2,q^3)$ are those arising from $\sigma$-fixed-planes of $\PG(8,q)$.

\begin{lemma}\Label{408} Let  ${\bar{\mathcal B}}$ be a $\phi$-fixed-$\Fq$-subplane of $\PG(2,q^3)$.
\begin{enumerate}
\item If $\g=1$, then there is a unique $\sigma$-fixed-plane $\beta$ in $\PG(8,q)$ with $\B(\beta)=\bar{\mathcal B}$.
\item If $\g=3$, then either
\begin{enumerate}
\item there are three $\sigma$-fixed-planes $\beta_1,\beta_2,\beta_3$ in $\PG(8,q)$ with $\B(\beta_i)=\bar{\mathcal B}$, $i=1,2,3$; or
\item $\bar{\mathcal B}$ contains only Type-\III-points and there are no $\sigma$-fixed-planes of $\PG(8,q)$ corresponding to $\bar{\mathcal B}$.
\end{enumerate}
\end{enumerate}
 \end{lemma}

\begin{proof} Let  ${\bar{\mathcal B}}$ be a $\phi$-fixed-$\Fq$-subplane of $\PG(2,q^3)$. So in $\PG(8,q)$, $\Bo{\mathcal B}$ is one ruling system of planes of a Segre variety $\S_{2;2}$. As $\Bo{\mathcal B}$  is fixed by $\sigma$, $\sigma$ permutes the $q^2+q+1$    $\Bo{{\mathcal B}}$-ruling planes.

Suppose there is a   $\sigma$-fixed $\Bo{{\mathcal B}}$-ruling plane $\beta$, so $\B(\beta)=\bar{\mathcal B}$. As $\bar{\mathcal B}$ is an $\Fq$-plane, we see from Table~\ref{tableB} that $\beta$ is either $\piq$, a fixed-\III-plane or an \hb-plane.  In particular, by Theorems~\ref{765} and \ref{113}, $\beta$ contains at least one fixed point, denote this by $P$. The  point $P$ lies in a unique $\SI$-plane $\alpha$ which is a plane of $\Bo{\mathcal B}$. By Corollary~\ref{294b}, $\alpha$ contains exactly $\g$ fixed points. Hence there are at most $\g$ $\sigma$-fixed $\Bo{{\mathcal B}}$-ruling planes.

Suppose $\g=1$, so $q^2+q+1\equiv 1\pmod3$. As $\sigma$ has order $3$,   $\sigma$ fixes at least one $\Bo{{\mathcal B}}$-ruling plane, and so fixes exactly one as required.

Suppose $\g=3$, so $q^2+q+1\equiv 3\pmod3$, so $\sigma$ fixes either zero or three $\Bo{{\mathcal B}}$-ruling planes.
Thus if $\sigma$ fixes a $\Bo{{\mathcal B}}$-ruling plane, then it fixes exactly three. As noted above, in this case, $\bar{\mathcal B}$ contains a Type-\I-point.  The remaining case is when $\sigma$ fixes zero $\Bo{{\mathcal B}}$-ruling planes and $\bar{\mathcal B}$ has no Type-\I-points.
Suppose $\bar{\mathcal B}$ contains no
 Type-\I-points. As $\phi$ fixes $\bar{\mathcal B}$, it induces a collineation acting on $\bar{\mathcal B}$ with no fixed points, hence no fixed lines. If $\bar{\mathcal B}$ contains a Type-\II-point $\bar B$, then $\bar B,\bar B^\phi,\bar B^{\phi^2}$ lie on a Type-\I-line and lie in $\bar{\mathcal B}$, contradicting $\bar{\mathcal B}$ having no fixed line.
 Hence $\bar{\mathcal B}$ contains only Type-\III-points as required.
   \end{proof}

   We note that $\Fq$-planes of $\PG(2,q^3)$ of the form described in Lemma~\ref{408}(2b) do exist, see \cite{BHJ2}.


\section{The Bruck-Bose $\PG(6,q)$ setting}\Label{Sec5}

We include a short discussion on the representation of the collineation $\phi$ in the Bruck-Bose $\PG(6,q)$ setting. The Bruck-Bose representation of $\PG(2,q^3)$ in $\PG(6,q)$ is discussed in detail in \cite{BJ-FFA}. We look at the case where the line at infinity $\li$ is a Type-\I-line. The analysis will be different if $\li$ has Type \II\ or \III.

Let $\Sigma_\infty$ be a hyperplane of $\PG(6,q)$. We call a subspace of $\PG(6,q)$ an \emph{affine subspace} if it is not contained in $\si$.
 Let $\mathcal S$ be a regular 2-spread of $\Sigma_\infty$. Let $\A(\mathcal S)$ be the incidence structure whose: points are the affine points of $\PG( 6,q)\setminus\si$; lines are the affine $3$-spaces which contain a plane of $\mathcal S$; with natural incidence. Then $\A(\mathcal S)\cong\textup{AG}(2,q^3)$.  We can uniquely complete this to a projective plane $\P(\mathcal S)\cong\PG(2,q^3)$, the points on the line at infinity $\li$ in $\P(\mathcal S)$ are in one to one correspondence with the planes of $\mathcal S$. If $\X$ is a subset of points of $\PG(2,q^3)$, then we denote the corresponding set of points in $\PG(6,q)$ by $[\X]$.  We can also construct this $\PG(6,q)$ representation by intersecting the $\PG(8,q)$ representation with
an appropriate $6$-space. Let $\Sigma_\infty$ be an $\H$-$5$-space of $\PG(8,q)$ and let $\Pi_6$ be a $6$-space containing $\Sigma_\infty$.
The regular 2-spread $\S$ meets $\si$ in a 2-spread $\mathcal S=\S\cap\si$.
Each $\S$-plane not contained in $\si$ meets $\PG(6,q)$ in a point, giving us the Bruck-Bose $\PG(6,q)$ representation.

 So if $\bar P$ is an affine point of $\PG(2,q^3)$, and $\bar P$ is a Type-$X$-point, $X\in\{\I,\II,\III\}$, then $[P]$ is an affine point of $\PG(6,q)$ and we call it an $X$-point.
If $\bar Q$ is a point of $\li$, then $[Q]$ is an $\mathcal S$-plane and we use the same naming as we did in $\PG(8,q)$. In particular, as $\li$ is a Type-\I-line, $\Sigma_\infty$ is  an  $\HI$-$5$-space.
The collineation $\phi$ of $\PG(2,q^3)$ corresponds to a projectivity  $\varphi$ of $\PG(6,q)$.

We can generalise  and use our analysis in $\PG(8,q)$ to  classify points,  lines and  planes of  $\PG(6,q)$ fixed by $\varphi$.  As the proofs are very similar to those in $\PG(8,q)$, we just   state the interesting classifications here, but omit the proofs.

First note that $\si$ corresponds to a Type-I-line, so has the same properties as  an  $\HI$-$5$-space in $\PG(8,q)$. Hence the $5$-space $\si$ contains exactly $\g(q+1)$ fixed points, $\g(q+1)(q-\m)$ $\I_\line$-points, $\g(q^3-q)$ $\II_\line$-points, $\g$ ptwise-fixed lines (each rules the $2$-regulus of $\SI$-planes in $\si$, $\g(q+1)$ fixed-\I-lines (which form $\g$ 1-reguli), and $\g(q^3-q)(q-\m)$ fixed-\II-lines. The fixed-$\I$-lines and fixed-\II-lines lie in a $3$-space,
 and form a linear congruence  as described in Theorem~\ref{191}.
It is then straightforward to classify the $\varphi$-fixed points in $\PG(6,q)$.
In particular, there are exactly $q^2+\g q+\g$ $\varphi$-fixed points in $\PG(6,q)$. These comprise the $q^2+q+1$ points in one $\varphi$-ptwise-fixed-plane $\pifix=[\piq]$, and an additional $\g-1$ $\varphi$-ptwise-fixed lines in $\si$.

A particularly interesting result is the following which looks at the set of all affine \I-points and affine \II-points, and shows that they are exactly the affine points of a quadric.

\begin{theorem}\Label{purplequadric} The set of all affine $\I$-points and $\II$-points in $\PG(6,q)$ is precisely the set of affine points on a quadric which is degenerate with vertex the plane $\pifix$ and base a hyperbolic quadric.
\end{theorem}

\begin{proof} We sketch the proof. Without loss of generality, let $\li$ have equation $z=0$,
let $\piq$ be the $\Fq$-plane $\piq =\{(x,y,z)\,|\, x,y,z\in\Fq, \textup{ not all }0\}$ of $\PG(2,q^3)$, and $$\varphi\colon (x,y,z)\mapsto (x^q,y^q,z^q).$$

The $\Fq$-plane $\piq$ corresponds to the plane $\pifix=[\piq]$ in $\PG(6,q)$, where $\pifix$ contains the  points
 with homogeneous coordinates $\{(a,0,0,b,0,0,c)\,|\, a,b,c\in\Fq\}.$
Let $\tau$ be a primitive element of $\Fqqq$ over $\Fq$. So each element $x\in\Fqqq$ can be uniquely written as $x=x_0+x_1\tau+x_2\tau^2$ with $x_0,x_1,x_2\in\Fq$.
Define the following two maps
$$\begin{array}{rrclcrrcl}
\theta:&\Fqqq&\longrightarrow &\Fq^3&\phantom{XXXXXXXX}&\bar\Theta:&\Fqqq\times\Fqqq\times\{1\}&\longrightarrow &\Fq^7\\
&x&\mapsto&(x_0,x_1,x_2)&
&&(x,y,1)&\mapsto&(\theta(x),\theta(y),1).
\end{array}
$$

An
 affine point $\bar P$ in $\PG(2,q^3)$ has homogeneous coordinates  $\vec P=(x,y,1)$. The corresponding affine point in the Bruck-Bose representation in  $\PG(6,q)$ has coordinates $$[P]=\bar\Theta( x, y,1) = (\theta(x),\theta(y),1).$$
The point $\bar P$
 has Type \I\ or Type \II\  iff $(x,y,1)$, $\phi(x,y,1)$, $\phi^2(x,y,1)$ are linearly dependent. That is,
\begin{eqnarray*} \label{purple}
\begin{vmatrix} x&y&1\\x^q&y^q&1\\
x^{q^2}&y^{q^2}&1
\end{vmatrix}
&=&0,
\end{eqnarray*}
and so,
\begin{eqnarray}\label{eqn-purple}
f(x,y)= (xy^q+x^qy^{q^2}+x^{q^2}y)-(xy^{q^2}+x^qy+x^{q^2}y^q)=0.
\end{eqnarray}
We write $x=x_0+x_1\tau +x_2\tau ^2$ and $y=y_0+y_1\tau +y_2\tau ^2$ for $x_i,y_i\in\Fq$.
  We first show that $f$ has no term in $x_0$ or $y_0$. From (\ref{eqn-purple}) the coefficient of $x_0$ is $(y^q+y^{q^2}+y)-(y^{q^2}+y+y^q)=0$. Similarly for $y_0$.

  Write $f(x,y)=\sum_{i=0}^2 (xy^q-x^qy)^{q^i}$. As $f(x,y)$ has no term in $x_0,y_0$, it follows that $xy^q-x^qy=(x_1\tau +x_2\tau ^2)(y_1\tau ^{q}+y_2\tau ^{2q})-(x_1\tau^q +x_2\tau ^{2q})(y_1\tau+y_2\tau ^{2})=(x_1y_2-x_2y_1)(\tau^{2q}\tau-\tau^2\tau^q)$. Thus we have $f(x,y)=(x_1y_2-x_2y_1)\times \beta$.
  As $f(x,y)=Tr(xy^q-x^qy)$, we have $f(x,y)\in\Fq$, so $\beta\in\Fq$.
   Further $\beta\ne 0$ as the \III-points do not lie on this equation.  Thus  this equation represents the quadric $\Q$ with equation $x_1y_2-x_2y_1$.
The singular space of $\Q$ is $x_1=x_2=y_1=y_2=0$, that is, the plane $\pifix$. To determine the type of base quadric, we intersect $\Q$ with the  $3$-space $\Sigma_3$ defined by $x_0=y_0=z=0$ which is disjoint from $\pifix$, the resulting quadric is non-singular hyperbolic quadric.
\end{proof}

By looking at how this quadric meets the $3$-spaces corresponding to a line of $\PG(2,q^3)$, we can determine the structure of the affine $\II$-points.

\begin{lemma}\Label{1133}
Let $\Q$ be the quadric whose affine points are precisely the affine $\I$-points and affine $\II$-points.
\begin{enumerate}
\item $\Q$ contains each $\SI$-plane of $\si$ and meets each $\SII$-plane of $\si$ in a non-degenerate conic.
\item Let $\Sigma_3$ be a $3$-space that meets $\si$ in an $\SI$-plane and corresponds to a Type-\II-line of $\PG(2,q^3)$. Then $\Sigma_3$ contains exactly $q^2$ affine $\II$-points which form an affine plane.
\item Let $\Sigma_3$ be a $3$-space that meets $\si$ in an $\SII$-plane and corresponds to a Type-\II-line of $\PG(2,q^3)$. Then $\Sigma_3$ contains exactly one affine $\I$-point $P$ and $q^2-1$ affine $\II$-points, these form the affine part of a quadratic cone with vertex $P$.
\item Let $\Sigma_3$ be a $3$-space that meets $\si$  in an $\SII$-plane and corresponds to a Type-\III-line of $\PG(2,q^3)$. Then $\Sigma_3$ contains exactly $q^2+q$ affine $\II$-points which form the affine part of a hyperbolic quadric.
\end{enumerate}
\end{lemma}

We conclude this section by briefly discuss the affine $\varphi$-fixed lines, there are two types, namely $\varphi$-ptwise-fixed lines contained in $\pifix$, and $\varphi$-fixed-\II-lines which exist iff $q\not\equiv-1\pmod3$. In this case, they occur in precisely the following way. Let $m$ be a $\varphi$-fixed-\I-line (so $m$ lies in an $\SI$-plane contained in $\si$) that meets $\pifix$ in a point. Note that $m$ exists iff $q\not\equiv -1\pmod3$.  Let $\ell$ be an affine line of $\pifix$ such that $m$ and $\ell$ meet. Then the lines through $m\cap\ell$ in the plane $\langle m,\ell\rangle$ distinct from $m,\ell$ are $\varphi$-fixed-$\II$-lines.

\section{Application to the Figueroa plane}\Label{Sec6}

In this section we define the Figueroa plane following the technique of \cite{grundhofer}, using the notation of \cite{BHJ2}.
As the Figueroa plane and the Desarguesian plane $\PG(2,q^3)$ have the same set of points, we can consider a  line of the Figueroa plane as a set of points of $\PG(2,q^3)$, and look at this set
in the $\PG(8,q)$ field reduction setting.

In $\PG(2,q^3)$, let $\PI$ denote the set of Type-\I-points; $\PII$ the set of Type-\II-points; and $\PIII$ the set of Type-\III-points. Similarly, let $\LI$ denote the set of Type-\I-lines; $\LII$ the set of Type-\II-lines; and $\LIII$ the set of Type-\III-lines.
If  we view lines as subsets of points, then $\PG(2,q^3)$ is an incidence structure   with points $ \PI\cup \PII\cup \PIII$  and lines $ \LI\cup \LII\cup \LIII$.  We define
the \emph{Figueroa plane} $\FIG(q^3)$, $q>2$, as the incidence structure   with points $\mathbb P=\PI\cup \PII\cup \PIII$  and lines $\LI\cup \LII\cup \LFig$, where $\LFig$ is defined below.
That is, lines of  $\LI\cup \LII$ contains the same set of points in both $\PG(2,q^3)$ and in $\FIG(q^3)$.
Elements of $\LFig$ are called
\emph{Fig-blocks}, and we view   a Fig-block as a subset of $q^3+1$  points in $\PG(2,q^3)$.
Each Type  $\III$ line in $ \LIII$ gives rise to a unique Fig-block in $\LFig$. Following the notation in \cite{BHJ2}, for a Type-\III-point $\bar G$, we let $\figg {\bar G}$ denote the Fig-block associated with the line $\bar G^\phi\bar G^{\phi^2}$, and $\LFig=\{\figg{\bar G}\,|\,\bar G\textup{ is a Type \III\ point}\}$. The set of points in the Fig-block $\figg{\bar G}$ is defined as $\figg {\bar G}=\bar \E_G\cup\bar \F_G$  where   \begin{itemize}
 \item  $\bar \E_G$ is the set of $q^2+q+1$ Type-\II-points on the line $\bar G^\phi \bar G^{\phi^2}$,
 \item $\bar \F_G=\{ \bar \ell^\phi\cap\bar \ell^{\phi^2}\,| \,\bar \ell \textup{ is a Type-\III-line through } \bar G\}$.
  \end{itemize}
 Note that the point $\bar G$ does not lie on the line $\bar G^\phi \bar G^{\phi^2}$ or in the Fig-block $\figg{\bar G}$; whereas  the
 points $\bar G^\phi, \bar G^{\phi^2}$ lie in both the line $\bar G^\phi \bar G^{\phi^2}$ and
in the Fig-block  $\figg{\bar G}$.

As we view   a Fig-block as a subset of $q^3+1$  points in $\PG(2,q^3)$, we can look at it in the $\PG(8,q)$ field reduction setting.
Under field reduction,   $\Bo{\E_G}$ is a set of $q^2+q+1$ $\SII$-planes in $\PG(8,q)$ lying in  an  $\HIII$-$5$-space. In fact, $\bar \E_G$ is a scattered linear set of rank $3$  and can be represented by a single plane of $\PG(8,q)$; Theorem~\ref{407} shows how there is a special way to construct such a plane using $\II_\line$-points.
The set $\bar \F_G$ contains the Type \III\ points in $\figg{\bar G}$, it corresponds under field reduction to a set $\Bo{ \F_G}$ of $q^3-q^2-q$ $\SIII$-planes in $\PG(8,q)$. In the case when $\g=\gcd(3,q-1)=1$,  that is $q\not\equiv 1\pmod3$,   we  construct a much smaller  geometrical set in $\PG(8,q)$ that corresponds to $\bar \F_G$.

We first generalise the notation used for linear sets in Section~\ref{Sec3}:  for \emph{any} pointset $\K$ in $\PG(8,q)$, we define a corresponding set $\B(\K)$ of points in $\PG(2,q^3)$ where $$\B(\K)=\{\bar P\in\PG(2,q^3)\,|\, \K\cap\Bo{P}\neq \emptyset\}.$$
A set ${\mathcal D}$ of $q+1$ planes contained in one ruling system of planes of a Segre variety $\S_{2;2}$ with the planes in ${\mathcal D}$ ruling non-degenerate conics is called a
  \emph{normal rational $3$-fold scroll} in   \cite[p93]{harris}.

\begin{figure}[h]
\hspace*{2cm}
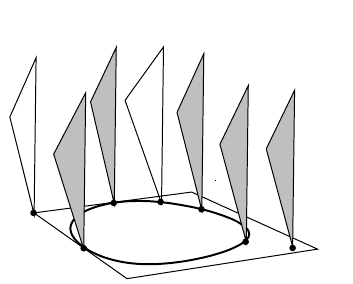\vspace*{5mm}
 \caption{A representation of the Fig-block $\figg{\bar G}$ in $\PG(8,q)$. By Theorem~\ref{407}, $\figg{\bar G}=\B(\beta\cup{\mathcal D}\setminus\pifix)$; the $q+1$ planes in $\beta\cup{\mathcal D}\setminus\pifix$ are shaded.}
 \label{fig4}
 \end{figure}

  The next result shows how we can represent a Fig-block of $\PG(2,q^3)$ using a normal rational $3$-fold scroll in $\PG(8,q)$. The construction
is illustrated in Figure~\ref{fig4}.

\begin{theorem}\Label{407} Suppose $q\not\equiv 1 \pmod3$. Let $\bar G$ be a Type-\III-point in $\PG(2,q^3)$ and $\bar m_G=\bar G^\phi \bar G^{\phi^2}$, and consider the $\PG(8,q)$ field reduction setting.
\begin{enumerate}
\item\label{407F} The set $\mathbb G$ of fixed-\III-planes meeting $\gamma=\Bo{G}$ forms one system of ruling planes of a Segre variety $\mathbb S_{2;2}$, denote the other system of ruling planes by $\mathbb G'$.
Let $P\in\gamma$ and $\pi=\langle P,P^\sigma,P^{\sigma^2}\rangle$, then $\pi\in\mathbb G$ and $\pi$ contains a unique fixed point $A$. Let $\mathcal C$ be the unique non-degenerate conic in $\pi$ that contains the three points $A,P^\sigma,P^{\sigma^2}$  and has   tangents $PP^\sigma$, $PP^{\sigma^2}$. The set $\D$ of $q+1$ planes of $\mathbb G'$ that meet $\C$ is a  normal rational 3-fold scroll  and  $\B({\mathcal D}\setminus\pifix)=\bar\F_G$.
\item\label{407E} The set of all $\II_\line$-points in the $\HIII$-$5$-space $\Bo{m_G}$ forms a plane $\beta\in\mathbb G'$
and   $\B(\beta)=\bar\E_G$.

\item $\figg{\bar G}= \B(\beta\cup{\mathcal D}\setminus\pifix)$ and  each point in $\figg{\bar G}$ distinct from $\bar G^\phi, \bar G^{\phi^2}$ corresponds to a unique point of $\beta\cup{\mathcal D}\setminus\pifix$.
\end{enumerate}
\end{theorem}

\begin{proof} By Corollary~\ref{333b1}, each of the $q^2+q+1$ $\SII$-planes in $\Sigma_5=\Bo{m_G}$ contains exactly one $\II_\line$-point, and so $\Sigma_5$ contains exactly $q^2+q+1$ $\II_\line$-points.
The condition $q\not\equiv 1\pmod{3}$ is equivalent to the condition $\g=1$, so  there is a unique hwise-fixed-$5$-space. By Corollary~\ref{221}, this hwise-fixed-$5$-space meets the $\HIII$-space $\Sigma_5$ in a plane $\beta$ that contains only $\II_\line$-points. That is, the $q^2+q+1$ $\II_\line$-points in $\Sigma_5$ form a plane $\beta$, and $\beta$ meets each $\SII$-plane in $\Sigma_5$ in a point.
Thus $\B(\beta)=\bar \E_G$.

As  $\bar G$ is a Type-\III-point, $\gamma=\Bo{G}$ is an $\SIII$-plane. By Theorem~\ref{401}, the set
 $\mathbb G$ of fixed-\III-planes meeting $\gamma$ is one system of ruling planes of a Segre variety $\S_{2;2}$, denote the other system of ruling planes by $\mathbb G'$.
Denote the planes of $\mathbb G$ by $\pi_i$, $i=1,\ldots,q^2+q+1$. By Theorem~\ref{765}, as $\g=1$, $\pi_i$ contains a unique fixed point, denote this by $A_i$. It follows from  Theorem~\ref{000}
that $\{A_i\,|\,i=1,\ldots,q^2+q+1\}$ is the unique pointwise fixed plane $\pifix$, so  $\pifix\in\mathbb G'$. Let $\pi_i\cap\gamma=P_i$ and let $\mathcal C_i$ be the unique non-degenerate conic in $\pi_i$ that contains the three points $A_i,P_i^\sigma,P_i^{\sigma^2}$  and has   tangents $P_iP_i^\sigma$, $P_iP_i^{\sigma^2}$. The conics $\C_i$, $i=1,\ldots,q^2+q+1$, are ruled by a set $\mathcal D$ of $q+1$ planes of $\mathbb G'$, moreover $\gamma,\gamma^\sigma,\gamma^{\sigma^2},\pifix\in\mathcal D$.

As  $\pi_i\notin\mathbb S$, the sets $\B(\pi_i)$, $i=1,\ldots,q^2+q+1$, are $\Fq$-planes of $\PG(2,q^3)$.
As $\gamma,\gamma^\sigma,\gamma^{\sigma^2}\in\mathbb G'$, we have $\bar G, \bar G^\phi, \bar G^{\phi^2}\in \B(\pi_i)$, so $\B(\pi_i)$, $i=1,\ldots,q^2+q+1$, are $\sigma$-fixed-planes.
As $q\not\equiv1\pmod{3}$, by Lemma~\ref{408}, the $\Fq$-planes $\B(\pi_i)$, $i=1,\ldots,q^2+q+1$, are distinct. As $\pifix\in\mathbb G'$, we can denote the Type-\I-points of $\PG(2,q^3)$ by $\bar T_i$, $i=1,\ldots,q^2+q+1$, so that $\bar T_i\in\B(\pi_i)$.
In particular $\B(\pi_i)$ is the unique $\Fq$-plane containing the points $\bar G, \bar G^\phi, \bar G^{\phi^2},\bar T_i$.
Let $\bar{\mathcal O}_i$ be the unique $\Fq$-conic in $\B(\pi_i)$ containing the points $\bar T_i,\bar G^\phi, \bar G^{\phi^2}$ whose tangents at $\bar G^\phi$ and $ \bar G^{\phi^2}$ meet at $\bar G$.
As $\g=1$, we can use \cite[Thm 1]{BHJ2} to conclude that  $\bar {\mathcal F}_G={\cup_{i=1}^{q^2+q+1}} \bar {\mathcal O}_i \setminus \bar T_i$. By construction, $\B(A_i)=\bar T_i$ and $\B(\mathcal C_i)=\bar{\mathcal O}_i$, hence $\B({\mathcal D}\setminus\pifix)=\bar\F_G$. Moreover, each point in $\bar\F_G$ distinct from $\bar G^\phi$, $\bar G^{\phi^2}$ corresponds to a unique point of ${\mathcal D}\setminus\pifix$.
It also follows from \cite[Thm 1]{BHJ2} that  there is a unique point of $\bar\E_G$ contained in $\B(\pi_i)$, $i=1,\ldots,q^2+q+1$, denote this point by $\bar L_i$. In $\PG(8,q)$, $\Bo{L_i}$ is an $\SII$-plane that meeets $\pi_i$ in a point $X_i$, $i=1,\ldots,q^2+q+1$. By Theorem~\ref{765}, $X_i$ is a $\II_\line$-point.  By Corollary~\ref{333b1}, $\Bo{L_i}$ contains a unique $\II_\line$-point, so (recalling the construction of $\beta$ from the first paragraph) $X_i\in\beta$.  That is,  $\beta$ meets each of $\pi_i$, $i=1,\ldots,q^2+q+1$, and so $\beta\in\mathbb G'$.
\end{proof}

{\bfseries Remark.\ } The underlying geometry is very different in the case when $q\equiv 1 \pmod3$, and this construction does not work.
The different geometrical structure when $q\equiv 1 \pmod3$ is highlighted in  Lemma~\ref{408}, which shows
 that in this case, the $q^2+q+1$  $\sigma$-fixed-planes in $\mathbb G$ correspond to only $\frac13(q^2+q+1)$ $\phi$-fixed-$\Fq$-planes of $\PG(2,q^3)$.
(In the notation of the proof of Theorem~\ref{407}, the set $\{\B(\pi_i)\,|\,i=1,\ldots,q^2+q+1\}$ has size
$\frac13(q^2+q+1)$.) Also see the discussion in \cite{BHJ2} as to why the conic construction in $\PG(2,q^3)$ does not generalise to the case $q\equiv 1 \pmod3$.

\begin{table}[h]
\caption{For each given subspace, this table states the total number of such subspaces in $\PG(8,q)$}\begin{center}
\begin{tabular}{rl}
\toprule
 subspace&number\\
\midrule
$\SI$-plane&$q^2+q+1$\\
$\SII$-plane&$(q^2+q+1)(q^3-q)$\\
$\SIII$-plane&$ q^3(q-1)^2(q+1) $\\
\midrule
fixed point&$\g (q^2+q+1)$\\
$\I_\line$-point&$\g (q-\m )(q^2+q+1)$\\
$\I_\plane$-point&$(q^2+q+1-\g (q-\m +1))(q^2+q+1)$\\
$\II_\line$-point&$\g (q^2+q+1)(q^3-q)$\\
$\II_\plane$-point&$(q^2+q+1)(q^3-q)(q^2+q+1-\g)$\\
$\III_\plane$-point&$q^3(q^2-1)(q^3-1)$\\
\midrule
ptwise-fixed line&$\g (q^2+q+1)$\\
fixed-\I-line &$\g (q^2+q+1)$\\
fixed-\II-line &$\g (q^3-q)(q^2+q+1)/(q-\m )$\\
\midrule
ptwise-fixed plane&$\g $\\

fixed-\IIa-plane  &$\g (q-2+\g-\m)(q+1)(q^2+q+1)$\\
fixed-\IIb-plane
 &$(q^3-q)(q^2+q+1)$\\
fixed-\III-plane&$q^3(q^2-1)(q-1)(q^2+q+1)/(q^2+q+1-\g (q-\m +1))$ \\
$\ha$-plane   &$\g(\m+1)(q^2+q+1)(q+1)$\\
$\hb$-plane&$\g(\m+1)(q^2+q+1)(q^2-1+\m)$\\
 \midrule
 hwise-fixed-$5$-space&$\g $\\
\bottomrule
\end{tabular}
\end{center}
\label{table16}
\end{table}

\newcommand{\qmo}{{\mathsf k}}
\begin{table}
\caption{For each given subspace  in $\PG(8,q)$, this table gives the number of points of each type in the subspace (where ${\mathsf v}=q^2+q+1$, $\mathsf k=q-\m+1$)}
\begin{tabular}{l|llllll}
\toprule
{subspace}&fixed-point&$\I_\line$-point&$\I_\plane$-point&$\II_\line$-point&$\II_\plane$-point&$\III_\plane$-point\\
\midrule
ptwise-fixed line&$q+1$&$0$&$0$&$0$&$0$&$0$\\
fixed-\I-line &$\m +1$&$q-\m $&$0$&$0$&$0$&$0$\\
fixed-\II-line& $\m +1$&$0$&$0$&$q-\m $&$0$&$0$\\
\midrule
$\SI$-plane&$\g $&$\g (q-\m )$&${\mathsf v}-\g \qmo$&$0$&$0$&$0$\\
$\SII$-plane&$0$&$0$&$0$&$\g $&${\mathsf v}-\g $&$0$\\
$\SIII$-plane&$0$&$0$&$0$&$0$&$0$&${\mathsf v}$\\
\midrule
ptwise-fixed-plane&${\mathsf v}$&$0$&$0$&$0$&$0$&$0$\\
fixed-\IIa-plane &$\g $&$q-\m $&$0$&$(\g -1)(q-\m )$&${\mathsf v}-\g \qmo$&$0$\\
 fixed-\IIb-plane&$\g $&$0$&$0$&$\g (q-\m )$&${\mathsf v}-\g \qmo$&$0$\\
fixed-\III-plane& $\g $&$0$&$0$&$\g (q-\m )$&$0$&${\mathsf v}-\g \qmo$\\
$\ha$-plane&$q+1+\m$&$q-\m$&$0$&$q^2-q$&$0$&$0$\\
$\hb$-plane&$q+1+\m$&$0$&$0$&$q^2-\m$&$0$&$0$\\
\midrule
hwise-fixed-$5$-space&$(\m +1){\mathsf v}$&$(q-\m ){\mathsf v}$&$0$&$(q^3-q){\mathsf v}$&$0$&$0$\\
\bottomrule
\end{tabular}
\label{table14}
\end{table}

\begin{table}[h!]
\caption{For each given subspace  in $\PG(8,q)$, this table gives the number of fixed lines of each type in the subspace}
\begin{center}
\begin{tabular}{l|lll}
\toprule
{subspace}&ptwise-fixed line&fixed-\I-line&fixed-\II-line\\
\midrule
$\SI$-plane&$0$&$\g $&$0$\\
$\SII$-plane&$0$&$0$&$0$\\
$\SIII$-plane&$0$&$0$&$0$\\
\midrule
ptwise-fixed-plane&$q^2+q+1$&$0$&$0$\\
fixed-\IIa-plane &$0$&$1$&$\g-1 $\\
fixed-\IIb-plane &$0$&$0$&$\g $\\
fixed-\III-plane&$0$&$0$&$\g $\\
$\ha$-plane &$1$&$1$&$q+\m-1$\\
$\hb$-plane&$1$&$0$&$q+\m$\\
\midrule
$\HI$-$5$-space&$\g$&$\g(q+1)$&$\g(q^3-q)/(q-\m)$\\
$\HII$-$5$-space&$0$&$\g$&$0$\\
$\HIII$-$5$-space&$0$&$0$&$0$\\
\midrule
hwise-fixed-$5$-space&$(\m +1)(q^2+q+1)$&$q^2+q+1$&$(q^3-q)(q^2+q+1)/(q-\m )$\\
\bottomrule
\end{tabular}
\end{center}
\label{table18}
\end{table}

{\bfseries Author details}

S.G. Barwick. School of Mathematical Sciences, University of Adelaide, Adelaide, 5005, Australia.
susan.barwick@adelaide.edu.au

A.M.W. Hui. Applied Mathematics Program, BNU-HKBU United International College, Zhuhai, China / School of Mathematical and Statistical Sciences, Clemson University, USA, huimanwa@gmail.com

W.-A. Jackson. School of Mathematical Sciences, University of Adelaide, Adelaide, 5005, Australia.
wen.jackson@adelaide.edu.au

A.M.W. Hui acknowledges the support of National Natural Science Foundation of China (Grant No. 12071041).

\end{document}

%% file: fig4.pdf_t
\begin{picture}(0,0)%
\includegraphics{fig4.pdf}%
\end{picture}%
\setlength{\unitlength}{4144sp}%
\begin{picture}(2675,2144)(936,-2281)
\put(1043,-1069){\makebox(0,0)[lb]{\smash{\fontsize{10}{12}\normalfont {\color[rgb]{0,0,0}$\gamma$}%
}}}
\put(1449,-530){\makebox(0,0)[lb]{\smash{\fontsize{10}{12}\normalfont {\color[rgb]{0,0,0}\XXX}%
}}}
\put(3008,-1281){\makebox(0,0)[lb]{\smash{\fontsize{10}{12}\normalfont {\color[rgb]{0,0,0}$\beta$}%
}}}
\put(951,-311){\makebox(0,0)[lb]{\smash{\fontsize{14}{16.8}\normalfont {\color[rgb]{0,0,0}\YYY}%
}}}
\put(3596,-1406){\makebox(0,0)[lb]{\smash{\fontsize{12}{14.4}\normalfont {\color[rgb]{0,0,0}\GGG}%
}}}
\put(1369,-2219){\makebox(0,0)[lb]{\smash{\fontsize{10}{12}\normalfont {\color[rgb]{0,0,0}$P^{\sigma^2}$}%
}}}
\put(2024,-1812){\makebox(0,0)[lb]{\smash{\fontsize{10}{12}\normalfont {\color[rgb]{0,0,0}$A$}%
}}}
\put(1679,-1822){\makebox(0,0)[lb]{\smash{\fontsize{10}{12}\normalfont {\color[rgb]{0,0,0}$P^\sigma$}%
}}}
\put(1920,-957){\makebox(0,0)[lb]{\smash{\fontsize{10}{12}\normalfont {\color[rgb]{0,0,0}$\pifix$}%
}}}
\put(1060,-1905){\makebox(0,0)[lb]{\smash{\fontsize{10}{12}\normalfont {\color[rgb]{0,0,0}$P$}%
}}}
\put(1645,-941){\makebox(0,0)[lb]{\smash{\fontsize{10}{12}\normalfont {\color[rgb]{0,0,0}$\gamma^{\sigma}$}%
}}}
\put(1353,-1346){\makebox(0,0)[lb]{\smash{\fontsize{10}{12}\normalfont {\color[rgb]{0,0,0}$\gamma^{\sigma^2}$}%
}}}
\end{picture}%